\documentclass[11pt]{amsart}
\usepackage{amsmath}
\usepackage{amssymb}
\usepackage{epsf}
\usepackage{psfrag}
\usepackage[all]{xy}
\usepackage{latexsym}
\usepackage{amscd}
\newtheorem{Theorem}{Theorem}[section]
\newtheorem{Lemma}[Theorem]{Lemma}
\newtheorem{cor}[Theorem]{Corollary}

\newtheorem{Proposition}[Theorem]{Proposition}
\theoremstyle{definition}
\newtheorem{defn}[Theorem]{Definition}
\newtheorem{conv}[Theorem]{Convention}
\newtheorem{rem}[Theorem]{Remark}
\newtheorem{exmp}[Theorem]{Example}

\begin{document}

\title[Foldings, graphs of groups and the membership problem]
{Foldings, graphs of groups and the membership problem}

\author[I.~Kapovich]{Ilya Kapovich} \address{Dept. of Mathematics,
  University of Illinois at Urbana-Champaign, 1409 West Green Street,
  Urbana, IL 61801, USA} \email{kapovich@math.uiuc.edu}

\author[R.~Weidmann]{Richard Weidmann} \address{
Fachbereich Mathematik,
Johann Wolfgang Goethe-Universit\"{a}t,
Robert Mayer-Strasse 6-8,
60325 Frankfurt (Main),
Germany}
\email{rweidman@math.uni-frankfurt.de}

\author[A.~G.~Myasnikov]{Alexei Myasnikov} \address{Department of
  Mathematics, City College of CUNY, New York, NY 10031, USA}
\email{alexeim@att.net}

\thanks{The first author acknowledges the support of the
U.S.-Israel Binational Science Foundation grant no. 1999298}



\date{\today}

\begin{abstract}
  We introduce a combinatorial version of
  Stallings-Bestvina-Feighn-Dunwoody folding sequences. We then show
  how they are useful in analyzing the solvability of the uniform subgroup
  membership problem for fundamental groups of graphs of groups.
  Applications include coherent right-angled Artin groups and
  coherent solvable groups.
\end{abstract}

\maketitle

\section{Introduction}

The idea of using foldings to study group actions on trees was
introduced by Stallings in a seminal paper~\cite{Sta83}, where he
applied foldings to investigate free groups. Free groups are
exactly those groups that admit free actions on simplicial trees.
Later Stallings~\cite{Sta91} offered a way to extend these ideas
to non-free actions of groups on graphs and trees.  Bestvina and
Feighn~\cite{BF91} gave a systematic treatment of Stallings'
approach in the context of graphs of groups and applied this
theory to prove a far-reaching generalization of Dunwoody's
accessibility result. Later Dunwoody~\cite{Du98} refined the
theory by introducing vertex morphism.  Dunwoody~\cite{Du99} used
foldings to construct a small unstable action on an $\mathbb
R$-tree. Some other applications of foldings in the graph of
groups context can be found in
\cite{RS,Sela97,Sela01,DD,Du97,Gui98,Gui00,Bow,BW,GuSc}.

In this paper we develop a combinatorial treatment of foldings geared
towards more computational questions. In particular we are interested
in the subgroup membership problem and in computing algorithmically
the induced splitting for a subgroup of the fundamental group of a
graph of groups.  Recall that a finitely generated group
\[G=\langle x_1,\dots, x_k \, | \, r_1, r_2,\dots \rangle\] is said to
have \emph{solvable membership problem} (or \emph{solvable uniform
  membership problem}) if there is an algorithm which, for any finite
family of words $u, w_1,\dots, w_n$ in $\{x_1,\dots, x_k\}^{\pm 1}$
decides whether or not the element of $G$ represented by $u$ belongs
to the subgroup of $G$ generated by the elements of $G$ corresponding
to $w_1,\dots, w_n$ (it is easy to see that this definition does not
depend on the choice of a finite generating set for $G$). Similarly,
if $H\le G$ is a specific subgroup, then $H$ is said to have
\emph{solvable membership problem in $G$} if there is an algorithm
deciding for any word $u$ in $\{x_1,\dots, x_k\}^{\pm 1}$ whether $u$
represents an element of $H$.

Amalgamated free products, HNN-extensions and more generally,
fundamental groups of graphs of groups play a very important role in
group theory.  However, till now there has been relatively little
understanding of how these fundamental constructions affect the
subgroup membership problem.  One of the first results in this
direction is due to Mihailova, who proved~\cite{Mi59,Mi68} that if $A$
and $B$ have solvable membership problem then so does their free
product $A\ast B$ (see also the subsequent work of Boydron~\cite{Bo}).
Mihailova~\cite{Mi58} also produced some important counter-examples
demonstrating the difficulty of the membership problem. Namely, she
proved that the direct product $G=F(a,b)\times F(x,y)$ of two free
groups of rank two contains a finitely generated subgroup $H$ with
unsolvable membership problem in~$G$. The group $F(a,b)\times F(x,y)$ can be
thought of as a double HNN-extension of $F(a,b)$:
\[
G=\langle F(a,b), x,y\, |, x^{-1}fx=f, y^{-1}fy=f \text{ for any }
f\in F(a,b)\rangle.
\]
It is well-known that a finitely generated free group has uniform
membership problem solvable in quadratic time in terms of
$|u|+|w_1|+\dots +|w_n|$. Thus even seemingly innocuous free
constructions have the potential of greatly affecting the complexity
of the membership problem. Another important example which to this
date is not at all understood is that of the mapping torus of a free
group automorphism.

Namely, let $G$ be a group and let $\phi:G\to G$ be an automorphism of
$G$. Then the HNN-extension of $G$ along $\phi$
\[
M_{\phi}:=\langle\, G, t \,|\, t^{-1}gt=\phi(g), \text{ for every }
g\in G\, \rangle =G\rtimes_{\phi} {\mathbb Z}.
\]
is called the \emph{mapping torus group} of $\phi$.

The case when $G$ is a free group, or more generally, a surface group,
is of particular importance in 3-dimensional topology. Yet, apart from
a few obvious observations, nothing is known about the solvability of
the membership problem for mapping tori of automorphisms of free
groups and surface groups.

A substantial amount of work on the membership problem for amalgamated
products and HNN-extension was done by
Bezverkhnii~\cite{Bez1,Bez2,Bez3,Bez4}. However, he did not use the
machinery of Bass-Serre theory of graphs of groups and groups acting
on trees. Consequently, all of his results have to rely on Britton's
lemma and the normal form theorem for amalgamated products, which
makes his proofs extremely technical and statements of most results
quite special.

Our goal is to present a more geometric and unified approach to this
topic which relies on Bass-Serre theory~\cite{Serre,Bass} as well as
on combinatorial foldings methods.  Because of our algorithmic goals,
when approximating an induced splitting for a subgroup of the
fundamental group of a graph of groups $\mathbb A$, we need to work
primarily at the level of quotient graphs of groups rather than at the
level of the Bass-Serre covering trees, as it is done in the
Stallings-Bestvina-Feighn-Dunwoody treatment of foldings.  We use
finite combinatorial objects called \emph{$\mathbb A$-graphs} (where
$\mathbb A$ is a given graph of groups) to provide such
approximations. For algorithmic reasons, $\mathbb
A$-graphs are labeled by elements and subgroups of the original
vertex groups of $\mathbb A$, rather than by some abstractly defined
groups and their homomorphisms into the original vertex groups.  

The
full list of conditions that guarantee that foldings of $\mathbb A$-graphs can be applied algorithmicaly and terminate yielding the induced splitting of an arbitrary finitely generated subgroup of $\pi_1(\mathbb A)$ turns out to be rather cumbersome (see
Definition~\ref{benign}, Theorem~\ref{construct} and
Theorem~\ref{member} below). The same is true for the conditions that guarantee that the membership problem is solvable (Theorem~\ref{member}). Instead we formulate a
corollary of the main results:

\begin{Theorem}\label{A}
  Let $\mathbb A$ be a finite graph of groups such that:

\begin{enumerate}
\item For every vertex $v$ of $A$ the vertex group $A_v$ is either
  locally quasiconvex word-hyperbolic or polycyclic-by-finite.
\item Every edge group of $\mathbb A$ is polycyclic-by-finite.

\end{enumerate}

Then for any vertex $v_0\in VA$ the uniform membership problem for
$G=\pi_1({\mathbb A}, v_0)$ is solvable. Moreover there is an
algorithm which, given a finite subset $S\subseteq G$, constructs the
induced splitting and a finite presentation for the subgroup
$U=\langle S\rangle\le G$.
\end{Theorem}

By the \emph{induced splitting} of $U$ in the theorem above we mean
the decomposition of $U$ as $U=\pi_1(\mathbb B,v_0)$ where $\mathbb B$
is the quotient graph of groups for the action of $U$ on the minimal
$U$-invariant subtree of the Bass-Serre covering tree
$X=\widetilde{({\mathbb A}, v_0)}$ that contains the base-vertex of
$X$.

The above theorem applies to a wide variety of situations. For
example, it is applicable to a finite graph of groups where all vertex
groups are virtually abelian or where all vertex groups are virtually
free and edge groups are virtually cyclic. In particular, the mapping
torus of an automorphism of a free abelian group of finite rank (or in
fact of any virtually polycyclic group) falls into this category. While
Theorem~\ref{A} does not say anything about the computational
complexity of the algorithm solving the membership problem, we believe
that in many specific cases this complexity can be analyzed and
estimated explicitly. For example, in the case when all vertex groups
are free and edge groups are cyclic, the folding algorithm provided by
Theorem~\ref{A} appears to have polynomial complexity. Indeed, Paul
Schupp~\cite{Sch} obtained more precise results with polynomial
complexity estimates for multiple HNN-extensions of free groups with
cyclic associated subgroups.

Not surprisingly, we also recover (see Corollary~\ref{Mih}) a
generalization of Mihalailova's theorem regarding the membership
problem for free products to graphs of groups with finite edge groups.

As an illustration of the usefulness of Theorem~\ref{A}, we apply it
to graph products and right-angled Artin groups.  Recall that if
$\Gamma$ is a finite simple graph with a group $G_v$ associated to
each vertex of $\Gamma$ then the \emph{graph product} group $G$ is
defined as the free product $\ast_{v\in V\Gamma} G_v$ modulo the
relations $[G_v,G_u]=1$ whenever $u$ and $v$ are adjacent vertices in
$\Gamma$. If each $G_v$ is an infinite cyclic group,
then $G$ is called a \emph{right-angled Artin group} or {\em graph
  group} and is denoted by $G(\Gamma)$.

\begin{cor}\label{product}
  Let $T$ be a finite tree such that for every vertex $v\in VT$ there
  is an associated finitely generated virtually abelian group $G_v$.
  Then the graph product group $G$ has solvable uniform membership
  problem. Moreover, there is an algorithm which, given a finite
  subset $S\subseteq G$, constructs a finite presentation for the
  subgroup $U=\langle S\rangle \le G$.
\end{cor}
\begin{proof}
  Note that for any groups $K,H$ we can write the direct product
  $H\times K$ as an amalgam:
\[
H\times K= H \ast_{H} (H\times K) \ast_{K} K.
\]

Let $v_1,\dots, v_n$ be the vertices of $T$.  Let $T'$ be the
barycentric subdivision of $T$. We give $T'$ the structure of a graph
of groups as follows. For each vertex $v_i$ of $T$ assign the vertex
group $T'_{v_i}:=G_{v_i}$.  For each barycenter $v$ of an edge
$[v_i,v_j]$ of $T$ assign the vertex group $T'_v:=G_{v_i}\times
G_{v_j}$. Also, for $e_i=[v_i,v]\in ET'$ and $e_j=[v_j,v]\in ET'$ put
$T'_{e_i}:=G_{v_i}$ and $T'_{e_j}:=G_{v_j}$. Finally, we define the
 boundary monomorphisms $T'_{e_i}\to T_{v}'$ and $T'_{e_i}\to T_{v_i}'$ to be
the inclusion map $G_{v_i}\to G_{v_i}\times G_{v_j}$ and the identity
map $G_{v_i}\to
G_{v_i}$ respectively. This defines a graph of groups $\mathbb T'$
where all vertex groups are finitely generated virtually abelian.
Moreover, we have $G\cong \pi_1({\mathbb T'}, T')$.

Corollary~\ref{product} now follows from Theorem~\ref{A}.
\end{proof}

Theorem~\ref{A} also applies to many right-angled Artin groups:

\begin{cor}\label{chordal}
  Let $G=G(\Gamma)$ be a coherent right-angled Artin group.  Then $G$
  has solvable uniform membership problem. Moreover, there is an
  algorithm which, given a finite subset $S\subseteq G$, constructs a
  finite presentation for the subgroup $U=\langle S\rangle \le G$.
\end{cor}

\begin{proof}

  Recall that a simple graph is called \emph{chordal} if it does not
  possess a chord-free simple circuit of length $\ge 4$, that is for
  every simple circuit of length $\ge 4$ there are two
  non-neighboring vertices in the circuit which are adjacent in the
  graph. For example, every tree is a chordal graph. Chordal graphs
  are of particular importance in the theory of right-angled Artin
  groups since by a result of Droms~\cite{Droms} a right-angled Artin group $G(\Gamma)$ is coherent if and only if $\Gamma$ is chordal.

  Let $G=G(\Gamma)$ be a coherent right-angled Artin group based on a
  finite graph $\Gamma$. Hence $\Gamma$ is chordal.

  By the result of Mihailova about free products mentioned above we
  may assume that $\Gamma$ is connected. We will think about the
  vertices of $\Gamma$ as the generators of $G$.

  Recall that a vertex $v$ of a simple graph is called \emph{simplicial}
  if any two vertices adjacent to $v$ are joined by an edge. It is a
  well-known graph-theoretic fact that every chordal graph has a simplicial vertex (see, for
  example, Lemma~5.3.16 in \cite{West}).

  In order to establish the corollary we need the following:

\noindent{\bf Claim.} Let $\Gamma$ be a finite connected chordal graph. Then there exists a tree
of free abelian groups $\mathbb T$ with $G=G(\Gamma)=\pi_1({\mathbb T}, T)$
such that for every free Abelian subgroup $A$ of $G$ that corresponds to a complete subgraph of $\Gamma$ there is a vertex group of $\mathbb T$
containing $A$.

We will prove the Claim by induction on the number of vertices in
$\Gamma$. When this number is $1$ or $2$, the statement is trivial.
Suppose $|V\Gamma|=n>2$ and the Claim has been verified for all graphs
with fewer than $n$ vertices.

Let $v$ be a simplicial vertex of $\Gamma$ and let $\Gamma_0$ be the
graph obtained from $\Gamma$ by removing $v$ and all edges adjacent to
$v$. Then $\Gamma_0$ is a chordal graph defining a right-angled Artin
group $G_0$ that is canonically embedded in $G$. Let $S$ be the set of
vertices of $\Gamma$ adjacent to $v$. Since $v$ is simplicial, the set
$S$ spans a complete subgraph of $\Gamma$ (and of $\Gamma_0$) and thus
defines a free abelian subgroup $A$ of $G_0$ and of $G$. By the
inductive hypothesis we may represent $G_0$ as $G_0=\pi_1({\mathbb
  T}_0, T_0)$ where $\mathbb T_0$ is a tree of free abelian groups
satisfying the requirements of the Claim for $\Gamma_0$. In
particular, there is a vertex $x$ of $T_0$ with vertex group $B$
such that $A\le B$.

We now enlarge $T_0$ to a tree $T$ by attaching an extra edge $e$ with
origin $x$. We define the vertex group for the new vertex $t(e)$ to be $A\times
\langle v\rangle$ and the edge group of $e$ to be $A$ (here $\langle
v\rangle$ is the infinite cyclic group $G_v$). The
boundary monomorphisms for $e$ are defined as the obvious inclusions. This
produces a tree of groups $\mathbb T$. By comparing the presentations
for $G_0=\pi_1({\mathbb T}_0, T_0)$ and for $G$ we see that
$G=\pi_1({\mathbb T}, T)$.

Moreover, every complete subgraph of $\Gamma$ is either contained in
$\Gamma_0$ or it is contained in the complete subgraph in $\Gamma$
spanned by $S$ and $v$. Hence $\mathbb T$ satisfies the requirements
of the Claim for $\Gamma$, and the Claim is established.

The statement of Corollary~\ref{chordal} now follows from
Theorem~\ref{A}.
\end{proof}

The simplest non-coherent right-angled Artin group is $F(a,b)\times
F(x,y)$. This group is based on an ``empty square'', that is a simple
circuit of length four, which is also the simplest example of a
non-chordal graph.  By Mihailova's theorem $F(a,b)\times F(x,y)$ has
unsolvable membership problem. Thus the statement of
Corollary~\ref{chordal} need not hold for non-coherent right-angled
Artin groups.

Another easy corollary of Theorem~\ref{A} is:

\begin{cor}
Let $G$ be a finitely generated coherent solvable group. Then $G$ has
solvable uniform membership problem.
\end{cor}
\begin{proof}
By a result of Groves~\cite{Grov} and Bieri-Strebel~\cite{BiSt} if $G$ is a finitely generated
coherent solvable group then either $G$ is polycyclic or $G$
is an ascending HNN-extension of a polycyclic group.
Hence $G$ has solvable uniform membership problem by Theorem~\ref{A}.
\end{proof}

The first author is grateful to Alexandr Kostochka, Derek Robinson and
Paul Schupp for helpful discussions. The authors also thank the
referee for a careful reading of the paper and for many helpful comments
and suggestions.

\section{Graphs of groups, subgroups and induced splittings}

We refer the reader to the book of Serre~\cite{Serre} as well as to
\cite{Ba,Bass,Co,SW} for detailed background information regarding
groups acting on trees and Bass-Serre theory.

\begin{conv} [Graph of groups notations]
  Following Serre, we say that a \emph{graph} $A$ consists of a vertex
  set $VA$, edge set $EA$, the inverse-edge function ${}^{-1}:EA\to
  EA$ and two edge endpoint functions $t:EA\to VA$, $o:EA\to VA$ with
  the following properties:

\noindent (1) The function ${}^{-1}$ is a fixed-point free involution
on $EA$;

\noindent (2) For any $e\in EA$ we have $o(e)=t(e^{-1})$.

We call $e^{-1}$ the \emph{inverse edge} of $e$. For $e\in EA$ we call
$o(e)$ the \emph{initial vertex} of $e$ and we call $t(e)$ the
\emph{terminal vertex} of $e$.

An edge-path in $A$ is \emph{reduced} if it does not contain a subpath
of the form $e, e^{-1}$, where $e\in EA$.

If $T$ is a tree and $v_0, v_1$ are vertices of $T$, we will denote by
$[v_0,v_1]_T$ the \emph{$T$-geodesic path} from $v_0$ to $v_1$, that
is the unique reduced edge-path from $v_0$ to $v_1$ in $T$.

A \emph{graph-of-groups} $\mathbb A$ consists of an underlying graph
$A$ together with the following data. For each vertex $v\in VA$ there
is an associated \emph{vertex group} $A_v$ and for each edge $e\in EA$
there is an associated \emph{edge group} $A_e$. Every edge $e\in EA$
comes equipped with two \emph{boundary monomorphisms}
${\alpha}_e:A_e\to A_{o(e)}$ and $\omega_e:A_e\to A_{t(e)}$ for all
$e\in EA$.  If $e^{-1}$ is the inverse edge of $e$ then we assume that
$A_{e^{-1}}=A_e$, $\alpha_{e^{-1}}=\omega_e$ and
$\omega_{e^{-1}}=\alpha_e$.
\end{conv}

\begin{defn}[$\mathbb A$-paths]
  Recall that in Bass-Serre theory if $\mathbb A$ is a graph of
  groups, then an $\mathbb A$-\emph{path of length $k\ge 0$ from $v\in
    VA$ to $v'\in VA$} is a sequence
\[
p=a_0,e_1,a_1,\dots, e_k,a_k
\]
where $k\ge 0$ is an integer, $e_1,\dots, e_k$ is an edge-path in $A$
from $v\in VA$ to $v'\in VA$, where $a_0\in A_{v}, a_k\in A_{v'}$ and
$a_i\in A_{t(e_i)}=A_{o(e_{i+1})}$ for $0<i<k$. We will call $k$ the
\emph{length} of $p$ and denote it by denote by $|p|$.
Note that we allow $k=|p|$ to be equal to zero, in which case $v=v'$
and $p=a_0\in A_v$.

If $p$ is an $\mathbb A$-path from $v$ to $v'$ and $q$ is an $\mathbb
A$-path from $v'$ to $v''$, then the \emph{concatenation $pq$ of $p$
  and $q$} is defined in the obvious way and is an $\mathbb A$-path
from $v$ to $v''$ of length $|p|+|q|$.
\end{defn}

\medskip The following notion plays an important role in
Bass-Serre theory.

\begin{defn}[Fundamental group of a graph of groups]
  Let $\mathbb A$ be a graph of groups.  Let $\sim$ be the equivalence
  relation on the set of all $\mathbb A$-paths generated (modulo
  concatenation) by:
\[
a,e, \omega_e(c),e^{-1},\bar a\ \sim\ a\alpha_e(c)\bar a,\ \text{ where } e\in EA, c\in
A_e\hbox{ and }a,\bar a\in A_{o(e)}.
\]

If $p$ is an $\mathbb A$-path, we will denote the $\sim$-equivalence
class of $p$ by $\overline p$. Note that if $p\sim p'$ then $p,p'$
have the same initial vertex and the same terminal vertex in $VA$.

Let $v_0\in VA$ be a vertex of $A$. We define the \emph{fundamental
  group} $\pi_1({\mathbb A},v_0)$ as the set of $\sim$-equivalence
classes of $\mathbb A$-paths from $v_0$ to $v_0$. It can be shown that
$G$ is in fact a group with multiplication corresponding to
concatenation of paths.
\end{defn}

Suppose that an $\mathbb A$-path $p$ has a subsequence of the form $a,e,
\omega_e(c),e^{-1},\bar a$. Replacing this subsequence in $p$ by
$a\alpha_e(c)\bar a$ produces an $\mathbb
A$-path $q$. In this situation we will say that $q$ is obtained from
$p$ by an {\em elementary reduction}. Note that $|q|=|p|-2$ and that
$p\sim q$.
If no elementary reductions are applicable to $p$, we say that $p$ is
$\mathbb A$-{\em reduced} (or just \emph{reduced}).

Any $\mathbb A$-path is equivalent to a reduced $\mathbb A$-path, and
such a reduced $\mathbb A$-path can be obtained by applying elementary
reductions as long as possible. The following proposition implies that
the reduced $\mathbb A$-path obtained in this way is almost unique.

\begin{Proposition}[Normal Form Theorem]\label{form}
  Let $\mathbb A$ be a graph of groups. Then:

\begin{enumerate}
\item If $a\in A_v, a\ne 1$ is a nontrivial vertex group element then
  the length zero path $a$ from $v$ to $v$ is not $\sim$-equivalent to
  the trivial path $1$ from $v$ to $v$.

\item Suppose $p=a_0,e_1,a_1,\dots, e_k, a_k$ is a reduced $\mathbb
  A$-path from $v$ to $v'$ with $k>0$. Then $p$ is not
  $\sim$-equivalent to a shorter path from $v$ to $v'$. Moreover, if
  $p$ is equivalent to a reduced $\mathbb A$-path $p'$ from $v$ to
  $v'$ then $p'$ has underlying edge-path $e_1, e_2,\dots, e_k$.

\item Suppose $T$ is a maximal subtree of $A$ and let $v_0\in VA$ be a
  vertex of $V$. Let $G=\pi_1({\mathbb A}, v_0)$. For $x,y\in VA$ we
  denote by $[x,y]_T$ the $T$-geodesic edge-path in $T$. Then $G$ is
  generated by the set $\overline S$ where
\[
S=\underset{e\in EA-ET}{\cup} [v_0,o(e)]_T \,e\, [t(e),v_0]_T \bigcup
\underset{v\in VA}{\cup} [v_0,v]_T \,A_v\, [v,v_0]_T
\]
\end{enumerate}
\end{Proposition}

We also need to recall the explicit construction of the Bass-Serre
universal covering tree for a graph of groups.

\begin{defn}[Bass-Serre covering tree]\label{defn:tree}
  Let $\mathbb A$ be a graph of groups with base-vertex $v_0\in VA$.
  We define an equivalence relation $\approx$ on the set of $\mathbb
  A$-paths originating at $v_0$ by saying that $p\approx p'$ if

\begin{enumerate}
\item $p$ and $p'$ are both $\mathbb A$-path from $v_0$ to $v$ for
  some $v\in VA$ and
\item $p\sim p'a$ for some $a\in A_v$.
\end{enumerate}

For a $\mathbb A$-path $p$ from $v_0$ to $v$, we shall denote the
$\approx$-equivalence class of $p$ by $\overline{p} A_v$.

We now define the Bass-Serre tree $\widetilde{({\mathbb A}, v_0)}$ as
follows. The vertices of $\widetilde{({\mathbb A}, v_0)}$ are
$\approx$-equivalence classes of $\mathbb A$-paths originating at
$v_0$. Thus each vertex of $\widetilde{({\mathbb A}, v_0)}$ has the
form $\overline p A_v$, where $p$ is an $\mathbb A$-path from $v_0$ to
a vertex $v\in VA$. (Hence we can in fact choose $p$ to be already
$\mathbb A$-reduced and such that the last group-element in $p$ is
equal to $1$.)

Two vertices $x, x'$ of $\widetilde{({\mathbb A}, v_0)}$ are connected
by an edge if and only if we can express $x,x'$ as $x=\overline p A_v,
x'=\overline{pae} A_{v'}$, where $p$ is an $\mathbb A$-path from $v_0$
to $v$ and where $a\in A_v$, $e\in EA$ with $o(e)=v, t(e)=v'$.

It follows from Proposition~\ref{form} that $\widetilde{({\mathbb A},
  v_0)}$ is indeed a tree. This tree has a natural base-vertex, namely
$x_0=\bar 1 A_{v_0}$ corresponding to the $\approx$-equivalence class
of the trivial path $1$ from $v_0$ to $v_0$.

Moreover, the group $G=\pi_1({\mathbb A},v_0)$ has a natural
simplicial action on $\widetilde{({\mathbb A}, v_0)}$ defined as
follows:

If $g=\overline{q}\in G$ (where $q$ is an $\mathbb A$-path from $v_0$
to $v_0$) and $u=\overline{p}A_v$ (where $p$ is an $\mathbb A$-path
from $v_0$ to $v\in VA$), then $g\cdot u:=\overline{qp}A_v$. It is not
hard to check that the action is well-defined on the set of vertices
of $\widetilde{({\mathbb A}, v_0)}$ and that it preserves the
adjacency relation. Thus $G$ in fact has a canonical simplicial action
without inversions on $\widetilde{({\mathbb A}, v_0)}$.
\end{defn}
It follows from Proposition~\ref{form} that if $p$ is an $\mathbb
A$-path from $v_0$ to $v$ then the map $A_v\to G$, $a\mapsto
\overline{p a p^{-1}}$ is an embedding. Moreover, in this case the
$G$-stabilizer of the vertex $\overline p A_v$ of
$\widetilde{({\mathbb A}, v_0)}$ is equal to the image of the above
map, that is to $\overline{pA_vp^{-1}}$. Similarly, the $G$-stabilizer
of an edge in $\widetilde{({\mathbb A}, v_0)}$ connecting $\overline p
A_v$ to $\overline{pae} A_{v'}$ is equal to $\overline p
(a\alpha_e(A_e) a^{-1}) \overline {p^{-1}}$.

The following well-known statement is the heart of
Bass-Serre theory and provides a duality between group actions on
trees and fundamental groups of graphs of groups.

\begin{Proposition}\label{prop:qt}
  Let $U$ be a group acting on a simplicial tree $Y$ without
  inversions. Then the graph $B=Y/U$ has a natural
  graph-of-groups structure $\mathbb B$ such that $U$ is canonically
  isomorphic to $\pi_1({\mathbb B}, v_0')$ and $Y$ is
  $U$-equivariantly isomorphic to the universal covering Bass-Serre
  tree of $\mathbb B$ (here $v_0'$ is the image of $v_0$ in $B$).
\end{Proposition}

\begin{rem}\label{exp}
  We want to remind the reader of the explicit construction of $\mathbb
  B$.  Let $T_1\subseteq Y$ and $T_2\subseteq Y$ be subtrees of $Y$
  such that the following hold:

\begin{enumerate}
\item $T_1\subseteq T_2$.
\item $T_1$ is the lift of a maximal subtree of $Y/U$ to $Y$.
\item $T_2$ is a fundamental domain for the action of $U$ on $Y$, i.e.
  $UT_2=Y$ and no two distinct edges of $T_2$ are $U$-equivalent.
\item Every vertex $v\in VT_2-VT_1$ is connected to a vertex of $T_1$ by
  a single edge.
\end{enumerate}

This clearly implies that no two vertices of $T_1$ are $U$-equivalent,
that $U(VT_1)=VY$ and that for every vertex of $v\in VT_2-VT_1$ there is
a unique vertex $x(v)\in VT_1$ which is $U$-equivalent to $v$.

For each vertex $v\in VT_2-VT_1$ choose an element $t_v\in U$ such that
$t_vv =x(v)$.  The graph of groups $\mathbb B$ is then defined as
follows.

\begin{enumerate}
\item The graph $B=Y/U$ is obtained from $T_2$ by identifying $v$ with
  $x(v)$ for each vertex $v\in VT_2-VT_1$. Thus we can assume that $T_1$
  is a subgraph of $B$ (in fact a spanning tree of $B$) and that
  $v_0'=v_0$.  Similarly, we assume that $EB=ET_2$. For any edge
  $e=[z,v]$ of $T_2$ with $z\in T_1$ and $v\in VT_2-VT_1$, we set
  $o_B(e)=z$ and $t_B(e)=x(v)$.
\item For each vertex $v\in VT_1$ we set $B_v:=Stab_U(v)$, where
  $Stab_U(v)$ is the $U$-stabilizer of $v\in X$.
\item For each edge $e=[z,v]\in ET_2$ we set $B_e:=Stab_U(e)$.
\item For each edge $e=[z,v]\in ET_1$ the boundary monomorphisms
  $\alpha_e^B: B_e\to B_z$ and $\omega_e^B:B_e\to B_v$ are defined as
  inclusions of $Stab_U(e)$ in $Stab_U(z)$ and $Stab_U(v)$
  accordingly.
\item Suppose $e=[z,v]$ is an edge of $T_2$ with $z\in T_1$, $v\in
  VT_2-VT_1$.  We set the boundary monomorphism $\alpha_e^B: B_e\to B_z$
  to be the inclusion of $Stab_U(e)$ in $Stab_U(z)$. We set the
  boundary monomorphism $\omega_e^B:B_e\to B_{x(v)}$ to be the map
  $g\mapsto t_v g t_v^{-1}, g\in B_{e}$.
\end{enumerate}
\end{rem}

\medskip
\begin{defn}[Induced splitting]
  Let $\mathbb A$ be a graph of groups with a base-vertex $v_0$.  Let
  $G=\pi_1({\mathbb A}, v_0)$ and let $X=\widetilde{({\mathbb A},
    v_0)}$ be the universal Bass-Serre covering tree of the based
  graph-of-groups $({\mathbb A}, v_0)$.  Thus $X$ has a base-vertex
  $x_0$ mapping to $v_0$ under the natural quotient map.

  Suppose $U\le G$ is a subgroup of $G$ and $Y\subset X$ is a
  $U$-invariant subtree containing $x_0$. Then the graph-of-groups splitting $\mathbb
  B$ of $U$ obtained as in Proposition~\ref{prop:qt} on the quotient
  graph $B=Y/U$ is said to be \emph{an induced splitting of $U\le G$
    with respect to $Y$ corresponding to the splitting
    $G=\pi_1({\mathbb A}, v_0)$}.
\end{defn}

If $U$ acts on $X$ without a global fixed point then there is a
preferred choice of a $U$-invariant subtree of $X$, namely the
smallest $U$-invariant subtree containing $x_0$, which will be denoted
$X_{U,x_0}$ (or by $X_{U}$, if no confusion is possible):
\[
X_{U,x_0}=X_U:=\cup_{u\in U} [x_0, ux_0]
\]

\medskip Notice that because of the explicit construction of $\mathbb
B$ each vertex group of $\mathbb B$ fixes a vertex of $X$ and hence is
conjugate to a subgroup of a vertex group of $\mathbb A$. Similarly,
edge groups of $\mathbb B$ are conjugate to subgroups of edge groups
of $\mathbb A$.  In practice, when talking about induced splittings,
we will often choose $Y$ to be $X_{U,x_0}$.

\section{$\mathbb A$-graphs}

In this section we introduce the combinatorial notion of an $\mathbb
A$-graph. These $\mathbb A$-graphs will approximate induced splittings
of subgroups of $\pi_1({\mathbb A}, v_0)$. In good situations, namely
when an $\mathbb A$-graph is ``folded'', an induced splitting can be
directly read off the $\mathbb A$-graph.

\begin{defn}[$\mathbb A$-graph]
  Let $\mathbb A$ be a graph of groups. \emph{An $\mathbb A$-graph}
  $\mathcal B$ consists of an underlying graph $B$ with the following
  additional data:

\begin{enumerate}

\item A graph-morphism $[\,.\, ]: B\to A$.

\item Each vertex $u\in VB$ has an associated group $B_u$, where
  $B_u\le A_{[u]}$.

\item To each edge $f\in EB$ there are two associated group elements
  $f_{\alpha}\in A_{[o(f)]}$ and $f_{\omega}\in A_{[t(f)]}$ such that
  $(f^{-1})_{\alpha}=(f_{\omega})^{-1}$ for all $f\in EB$.
\end{enumerate}
\end{defn}

\begin{conv}
  If $f\in EB$ and $u\in VB$, we shall refer to $e=[f]\in EA$ and
  $v=[u]\in VA$ as \emph{ the type of $f$ and $u$} accordingly.  Also,
  especially when representing $\mathbb A$-graphs by pictures, we will
  sometimes say that an edge $f$ of an $\mathbb A$-graph $\mathcal B$
  has \emph{label} $(f_{\alpha},[f], f_{\omega})$. Similarly, we will
  say that a vertex $u\in VB$ has \emph{label} $(B_u,[u])$.
\end{conv}

\begin{figure}[here]
  \centerline{ \footnotesize \setlength{\unitlength}{.8cm}
\begin{picture}(16,4)
\put(12.5,2){\circle*{.07}}
\bezier{200}(12.5,2)(11.4,3)(10.6,3)
\bezier{200}(12.5,2)(11.4,1)(10.6,1)
\bezier{200}(10.6,3)(9.7,3)(9.7,2)
\bezier{200}(10.6,1)(9.7,1)(9.7,2)
\put(12.5,2){\line(4,3){2}}
\put(12.5,2){\line(4,-3){2}}
\put(14.5,3.5){\circle*{.07}}
\put(14.5,.5){\circle*{.07}}
\put(3.5,2){\circle*{.07}}
\bezier{200}(3.5,2)(2.4,3)(1.6,3)
\bezier{200}(3.5,2)(2.4,1)(1.6,1)
\bezier{200}(1.6,3)(.7,3)(.7,2)
\bezier{200}(1.6,1)(.7,1)(.7,2)
\put(3.5,2){\line(4,3){2}}
\put(3.5,2){\line(4,-3){2}}
\put(5.5,3.5){\circle*{.07}}
\put(5.5,.5){\circle*{.07}}
\put(4.5,2.75){\vector(4,3){.01}}
\put(4.5,1.25){\vector(4,-3){.01}}
\put(13.5,2.75){\vector(4,3){.01}}
\put(13.5,1.25){\vector(-4,3){.01}}
\put(.7,2){\vector(0,1){.01}}
\put(9.7,2){\vector(0,-1){.01}}
\put(5.1,3.6){$(B_3,v_3)$}
\put(14.1,3.6){$(B_3,v_3)$}
\put(5,.1){$(B_2,v_2)$}
\put(14,.1){$(B_2,v_2)$}
\put(3.8,1.88){$(B_1,v_1)$}
\put(12.8,1.88){$(B_1,v_1)$}
\put(3.2,3){$(a',e',b')$}
\put(12.2,3){$(a',e',b')$}
\put(3.3,.7){$(a,e,b)$}
\put(12.2,.7){$(b^{-1},e^{-1},a^{-1})$}
\put(-.8,1.75){$(\overline a,\overline e,\overline b)$}
\put(7.1,1.75){$(\overline{b}^{-1},\overline{e}^{-1},\overline{a}^{-1})$}
\end{picture}}
\caption{Two distinct diagrams associated to the same $\mathbb A$-graph}\label{distinct diagrams}
\end{figure}
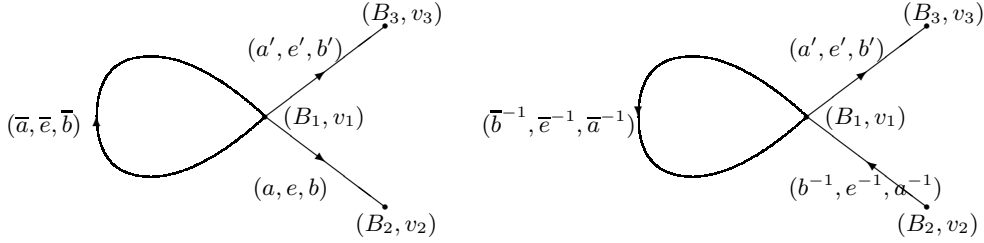

We will visualize an $\mathbb A$-graph $\mathcal B$ in the obvious way
by drawing the underlying graph $B$ with the appropriate labels next
to its vertices and edges. For every geometric edge we choose the
label of either edge of the corresponding edge-pair $\{f,f^{-1}\}$.
For convenience we will further orient the edge by attaching an arrow
such that for an edge with label $(a,e,b)$ one travels from a vertex
with label $(B,o(e))$ to a vertex with label $(B',t(e))$ if one
follows the direction of the arrow. It follows that reversing the
orientation of an edge and replacing the label $(a,e,b)$ by
$(b^{-1},e^{-1},a^{-1})$ yields another diagram of the same $\mathbb
A$-graph. An example is shown in Figure~\ref{distinct diagrams}.

To any $\mathbb A$-graph we can then associate in a natural way a
graph of groups:

\begin{defn}[Graph of groups defined by an $\mathbb A$-graph]

  Let $\mathcal B$ be an $\mathbb A$-graph. The associated graph of
  groups ${\mathbb B}$ is defined as follows:

\begin{enumerate}
\item The underlying graph of $\mathbb B$ is the graph $B$.

\item For each $u\in VB$ we put the vertex group of $u$ to be $B_u$.

\item For each $f\in EB$ we define the edge group of $f$ in $\mathbb
  B$ as
  $$B_f:=\alpha_{[f]}^{-1}(f_{\alpha}^{-1}B_{o(f)}f_{\alpha})\cap
  \omega_{[f]}^{-1}(f_{\omega}B_{t(f)}f_{\omega}^{-1})\le A_{[f]}.$$
\item For each $f\in EB$ we define the boundary monomorphism $\alpha_f:B_f\to
  B_{o(f)}$ as
  $\alpha_f(g)=f_{\alpha}\big(\alpha_{[f]}(g)\big)f_{\alpha}^{-1}$.
\end{enumerate}

\end{defn}

\begin{exmp}\label{edge}
  Let $\mathbb A$ be the ``edge-of-groups'' corresponding to an
  amalgamated product $G= P\ast_C Q$. Thus $A$ consists of an edge $e$
  with two distinct endpoints $v_0=o(e)$ and $v_1=t(e)$. The vertex
  and edge groups are: $A_{v_0}=P$, $A_{v_1}=Q$, $A_{e}=A_{e^{-1}}=C$.
  The boundary monomorphisms are the inclusions of $C$ into $P$ and $Q$.

  Consider the $\mathbb A$-graph $\mathcal B$, shown in
  Figure~\ref{Fi:edge}, consisting of a single edge $f$ of type $e$
  with $o(f)=u_0$ of type $v_0$ and $t(f)=u_1$ of type $v_1$. The
  associated group of $u_0$ is $P_0\le P$ and the associated group of
  $u_1$ is $Q_0\le Q$. Let $a=f_{\alpha}\in P$ and $b=f_{\omega}\in
  Q$. Thus the label of $f$ is $(a,e,b)$.

  Then the graph of groups $\mathbb B$ defined by $\mathcal B$ looks
  as follows. The underlying graph of $\mathbb B$ is still the single
  edge $f$. The vertex group of $u_0$ is $P_0$ and the vertex group of
  $u_1$ is $Q_0$. The edge group of $f$ is $B_f=a^{-1} P_0a\cap
  bQ_0b^{-1}\le C$. The boundary monomorphisms corresponding to $f$ are:
  $\alpha_f(c):=aca^{-1}$ and $\omega_f(c):=b^{-1}cb$ for $c\in B_f$.

\begin{figure}[here]
  \centerline{ \footnotesize \setlength{\unitlength}{.8cm}
\begin{picture}(16,2)
\put(2,1){\circle*{.07}}
\put(6,1){\circle*{.07}}
\put(10,1){\circle*{.07}}
\put(14,1){\circle*{.07}}
\put(2,1){\line(1,0){4}}
\put(10,1){\line(1,0){4}}
\put(4,1){\vector(1,0){.01}}
\put(12,1){\vector(1,0){.01}}
\put(1.5,1){$v_0$}
\put(1.1,.2){$A_{v_0}=P$}
\put(3.9,1.1){$e$}
\put(6.2,1){$v_1$}
\put(5.1,.2){$A_{v_1}=Q$}
\put(8.5,.9){$(P_0,v_0)$}
\put(11.3,1.2){$(a,e,b)$}
\put(14.2,.9){$(Q_0,v_1)$}
\end{picture}}
\caption{Example of an amalgamated product and an $\mathbb A$-graph}\label{Fi:edge}
\end{figure}

\end{exmp}

\smallskip\begin{conv} 

Let $\mathcal B$ be an $\mathbb A$-graph defining a graph-of-groups
$\mathbb B$.
Suppose $u,u'\in VB$ and $p$ is a $\mathbb
  B$-path from $u$ to $u'$. Thus $p$ has the form:
\[
p=b_0,f_1,b_1,\dots, f_s,b_s
\]
where $s\ge 0$ is an integer, $f_1,\dots, f_s$ is an edge path in $B$
from $u$ to $u'$, where $b_0\in B_u$, $b_s\in B_{u'}$ and $b_i\in
B_{t(f_i)}=B_{o(f_{i+1})}$ for $0<i<s$.  Recall that each edge $f_i$
has a label $(g_i,e_i,k_i)$ in $\mathcal B$, where $e_i=[f_i]$,
$g_i=(f_i)_{\alpha}$ and $k_i=(f_i)_{\omega}$.

Hence the $\mathbb B$-path $p$ determines the $\mathbb A$-path
$\mu(p)$ from $[u]$ to $[u']$ in $\mathbb A$ defined as follows:
\[
\mu(p)=(b_0g_1), e_1, (k_1b_1g_2), e_2,\dots, (k_{s-1}b_{s-1}g_s),e_s,
(k_sb_s)
\]
Notice that $|p|=|\mu(p)|$.
\end{conv}

We also want to think about an $\mathbb A$-graph as an ``automaton''
over $\mathbb A$ which ``accepts'' a certain subgroup of the
fundamental group of $\mathbb A$.
\begin{defn}
  Let $\mathcal B$ be an $\mathbb A$-graph with a base-vertex $u_0\in
  VB$.

  We define the \emph{language} $L({\mathcal B}, u_0)$ as
\[
L({\mathcal B}, u_0):=\{ \mu(p)| \, p\text{ is a reduced } {\mathbb
  B}\text{-path from } u_0 \text{ to } u_0 \text{ in } {\mathbb B}\}
\]
Thus $L({\mathcal B}, u_0)$ consists of $\mathbb A$-paths from
$v_0:=[u_0]$ to $v_0$.
\end{defn}

A simple but valuable observation states that the language of an
$\mathbb A$-graph represents a subgroup in the fundamental group of
$\mathbb A$.

\begin{Proposition}\label{basic}
  Let $\mathcal B$ be an $\mathbb A$-graph, $u_0\in VB$, $v_0=[u_0]$
  and $G=\pi_1({\mathbb A}, v_0)$.

  Then:

\noindent (1) If $p,p'$ are $\sim$-equivalent $\mathbb B$-paths, then
$\mu(p)\sim \mu(p')$ as $\mathbb A$-paths.

\noindent (2) The map $\mu$ restricted to the set of $\mathbb
B$-paths from $u_0$ to $u_0$ factors through to a homomorphism $\nu:
\pi_1({\mathbb B}, u_0)\to G$.

\noindent (3) We have
$\overline {L({\mathcal B}, u_0)}=\nu(\pi_1({\mathbb B}, u_0))$.  In
particular, $\overline {L({\mathcal B}, u_0)}$ is a subgroup of $G$.

\noindent (4) There is a canonical $\nu$-equivariant simplicial map
$\phi: \widetilde{({\mathbb B}, u_0)}\to \widetilde{({\mathbb A},
  v_0)}$ respecting the base-points.
\end{Proposition}

\begin{proof}
  Part (1) follows directly from the definitions of $\sim$ and
  $\mathbb B$.  Part (1) immediately implies parts (2) and (3).
  
  To establish (4) we will provide a direct construction of $\nu$
  which relies on the explicit definition of the Bass-Serre tree for a
  graph of groups given earlier. Denote $X=\widetilde{({\mathbb A},
    v_0)}$ and $Y=\widetilde{({\mathbb B}, u_0)}$. Let
  $y=\overline{p}B_u$ be a vertex of $Y$, where $p$ is a $\mathbb
  B$-path from $u_0$ to $u\in VB$. Denote $v=[u]\in VA$. We put
  $\phi(y):=\overline {\mu(p)}A_v\in VX$. First, note that this
  definition does not depend on the choice of $p$. Indeed, suppose
  $p'$ is another $\mathbb B$-path from $u_0$ to $u$ such that
  $p\approx p'$. Then by Definition~\ref{defn:tree}
  $\overline{p'}=\overline{pb}$ for some $b\in B_u\le A_v$. Hence
  $\overline{\mu(p)}A_v=\overline{\mu(p)}bA_v=\overline{\mu(pb)}A_v=\overline{\mu(p')}A_v$.
  Thus $\phi$ is well-defined on the vertex set of $Y$.
  
  It remains to check that $\phi$ preserves the adjacency relation.
  Let $y=\overline{p}B_u\in VY$ be as above and let
  $y'=\overline{pbf}B_{u'}\in VY$ be an adjacent vertex of $Y$, where
  $b\in B_u\le A_v$ and where $f\in EB$ is an edge of type $e\in EA$
  with $o(f)=u$. Thus $o(e)=v\in VA$. We already know that
  $\phi(y)=\overline{\mu(p)}A_v$. Denote $u'=t(f)$ and $v'=t(e)$, so
  that $[u']=v'$. Also denote $g=f_{\alpha}\in A_v$ and
  $h=f_{\omega}\in A_{v'}$.  Then $pbf$ is a $\mathbb B$-path from
  $u_0$ to $u'$.

  Therefore
\[
\phi(y')=\overline{\mu(pbf)}A_{v'}=\overline{\mu(p)bgeh}A_{v'}=\overline{\mu(p)bge}A_{v'}
\]
is an adjacent vertex of $\phi(y)=\overline{\mu(p)}A_v$ since $bg\in
A_v$.  Thus indeed $\phi$ is a well-defined simplicial map from $Y$ to
$X$.  We leave checking the equivariance properties of $\phi$ to the
reader.
\end{proof}

We will see that every subgroup $H$ of $G=\pi_1({\mathbb A}, v_0)$ arises
in this fashion, i.e. for every $H\le G$ we have $H=\nu(\pi_1(\mathbb B,u_0))$ where $\mathbb B$ is the graph of groups associated to some $\mathbb A$-graph $\mathcal B$.  Moreover, for an ``efficient'' choice of $\mathcal
B$ the homomorphism $\nu:\pi_1(\mathbb B,u_0)\to H$ is an isomorphism and the graph of groups $\mathbb B$ represents the
induced splitting of the subgroup $H\le G$ with respect to the action of $H$ on the Bass-Serre
covering tree of $\mathbb A$.

\begin{rem}
  Let $\mathbb A$ and $\mathcal B$ be as in Example~\ref{edge}. Then
\[
\overline {L({\mathcal B}, u_0)}=\nu(\pi_1({\mathbb B}, u_0))=\langle
P_0, abQ_0b^{-1}a^{-1}\rangle \le G=P \ast_C Q.
\]
\end{rem}

\medskip The following lemma is an immediate corollary of
Proposition~\ref{form} and Proposition~\ref{basic}:

\begin{Lemma}\label{gen}
  Let $\mathcal B$ be an $\mathbb A$-graph with a base-vertex $u_0$ of
  type $v_0$. Let $T\subseteq B$ be a spanning tree. For any two
  vertices $u,u'\in T$ denote by $[u,u']_T$ the $T$-geodesic path from
  $u$ to $u'$.

  Then $\pi_1({\mathbb B}, u_0)$ is generated by $\overline{S_T}$
  where $S_T$ is the following set:
\[
S_T:=\{[u_0,u]_T B_u [u,u_0]_T\ |\ u\in VB\}\cup \{[u_0,o(e)]_T e
[t(e),u_0]_T\ |\ e\in E(B-T)\}.
\]
In particular, $\overline {L({\mathcal B}, u_0)}\le\pi_1({\mathbb A},
v_0)$ is generated by $\overline{\mu(S_T)}=\nu(\overline{S_T})$.
\end{Lemma}

\section{Folding moves and folded $\mathbb A$-graphs}

\begin{defn}[Folded $\mathbb A$-graph]\label{defn:folded}
  Let $\mathcal B$ be an $\mathbb A$-graph.

  We will say that $\mathcal B$ is \emph{not folded} if at least one
  of the following applies:

\begin{enumerate}

\item There are two distinct edges $f_1,f_2$ with $o(f_1)=o(f_2)=z$
  and labels $(a_1,e,b_1)$, $(a_2,e,b_2)$ accordingly, such that $z$
  has label $(A',u)$ and $a_2=a'a_1\alpha_e(c)$ for some $c\in A_e$
  and $a'\in A'$.
\item There is an edge $f$ with label $(a,e,b)$, with $o(f)$ labeled
  $(A',u)$ and $t(f)$ labeled $(B',v)$ such that
  $\alpha_e^{-1}(a^{-1}A'a)\ne \omega_e^{-1}(bB'b^{-1})$.
\end{enumerate}

Otherwise we will say that $\mathcal B$ is \emph{folded}.
\end{defn}

It is easy to see that if $\mathcal B$ is folded then any reduced
$\mathbb B$-path translates into a reduced $\mathbb A$-path.

\begin{Lemma}
  Let $\mathcal B$ be a folded $\mathbb A$-graph defining the graph of
  groups $\mathbb B$. Suppose $p$ is a reduced $\mathbb B$-path.  Then
  the corresponding $\mathbb A$-path $\mu(p)$ is $\mathbb A$-reduced.
\end{Lemma}

\begin{proof}

  Suppose $p$ is a $\mathbb B$-reduced $\mathbb B$-path and $\mu(p)$
  is the corresponding $\mathbb A$-path. Assume that $\mu(p)$ is not
  reduced. Then $p$ has a subsequence of the form $f, a_1, f'$ where
  $f^{-1}, f'$ are edges of $B$ of the same type $e\in EA$ such that
  the label of $f^{-1}$ is $aeb$, the label of $f'$ is $a'eb'$, where
  $v\in VA$ is the type of $o(f')=t(f)\in VB$, $a,a'\in A_v$, $a_1\in
  B_{t(f)}\le A_v$ and the $\mathbb A$-path $b^{-1},e^{-1}, a^{-1}a_1a', e,b'$
  is not $\mathbb A$-reduced.
  
  This means that for some $c\in A_e$ we have
  $a^{-1}a_1a'=\alpha_e(c)$, that is $a'=a_1^{-1} a\alpha_e(c)$. If
  $f^{-1}$ and $f'$ are two distinct edges of $B$, this contradicts
  our assumption that $\mathcal B$ is folded. Thus $f^{-1}=f'$, so
  that $a=a', b=b'$. Therefore $a^{-1}a_1a=\alpha_e(c)$. Recall that
  since $\mathcal B$ is folded, part (2) of
  Definition~\ref{defn:folded} does not apply.  Therefore the edge
  group in $\mathbb B$ is $B_{f'}=\alpha_e^{-1}(a^{-1}A_1a)$, where
  $A_1=B_{t(f)}=B_{o(f')}$ and so $c\in B_{f'}$.  Moreover, the
   boundary monomorphism of $f'$ in $\mathbb B$ was defined as
  $\alpha_{f'}^B(c)=a\alpha_e(c)a^{-1}$. Thus
  $a_1=a\alpha_e(c)a^{-1}\in \alpha_{f'}^B(B_{f'})$.  Hence $f, a_1,
  f'$ is not $\mathbb B$-reduced, contrary to our assumptions.
\end{proof}

The above lemma immediately implies the following important fact:

\begin{Proposition}\label{prop:induced}
  Let $\mathcal B$ be a folded $\mathbb A$-graph defining the graph of
  groups $\mathbb B$. Let $u_0$ be a vertex of $B$ of type $v_0\in
  VA$. Denote $G=\pi_1({\mathbb A}, v_0)$ and $U=\overline
  {L({\mathcal B}, u_0)}=\nu(\pi_1({\mathbb B}, u_0))\le G$.

  Then the epimorphism $\nu: \pi_1({\mathbb B}, u_0) \to U$ is an
  isomorphism and the graph map $\phi$ between the Bass-Serre covering
  trees $\phi: \widetilde{({\mathbb B}, u_0)}\to \widetilde{({\mathbb
      A}, v_0)}$ is injective.
\end{Proposition}

Proposition~\ref{prop:induced} essentially says that if $\mathcal B$
is a folded $\mathbb A$-graph defining a subgroup $U\le G$, then
$U=\pi_1({\mathbb B}, u_0)$ is an induced splitting for $U\le
G=\pi_1({\mathbb A}, v_0)$.

\begin{exmp}

  Let $\mathbb A$ and $\mathcal B$ be as in Example~\ref{edge}. Recall
  that in this case $G=P\ast_C Q$ and $U= \langle P_0,
  abQ_0b^{-1}a^{-1}\rangle$. Recall also that in the graph of groups
  $\mathbb B$ the edge group of $f$ is $B_f=a^{-1}P_0a\cap
  bQ_0b^{-1}\le C$.

  By Definition~\ref{defn:folded} the $\mathbb A$-graph $\mathcal B$
  is folded if and only if $a^{-1}P_0a\cap C= bQ_0b^{-1}\cap C$ (in
  which case this last group is also equal to $B_f$). It is easy to
  see that, as claimed by Proposition~\ref{prop:induced}, if $\mathcal
  B$ is folded then $U=P_0 \ast_{aB_fa^{-1}} abQ_0b^{-1}a^{-1}$.
\end{exmp}

We will now describe certain moves, called \emph{folding moves} on
$\mathbb A$-graphs, which preserve the corresponding subgroups of the
fundamental group of $\mathbb A$. These folding moves are a more
combinatorial version of the folding moves of
Bestvina-Feighn~\cite{BF91} and Dunwoody~\cite{Du98}; implicitly they also contain Dunwoody's vertex morphisms.

\smallskip Whenever we make changes to the label of an edge $f$ of an
$\mathbb A$-graph we assume that the corresponding changes are made to
the label of $f^{-1}$.

\subsection{Auxiliary moves}

We will introduce three moves that can be applied to $\mathbb
A$-graphs. These moves do not substantially change its structure and
can be applied to any $\mathbb A$-graph.

\begin{defn}[Conjugation move $A0$]
  Let $\mathcal B$ be an $\mathbb A$-graph. Suppose $u$ is a vertex of
  $\mathcal B$ and that $g\in A_{[u]}$.

  Let $\mathcal B'$ be the $\mathbb A$-graph obtained from $\mathcal
  B$ as follows:

\begin{enumerate}
\item Replace $B_u$ by $gB_ug^{-1}$.
\item For each non-loop edge $f$ with $o(f)=u$ replace $f_{\alpha}$
  with $gf_{\alpha}$.
\item For each non-loop edge $t(f)=u$ replace $f_{\omega}$ with
  $f_{\omega}g^{-1}$.
\item For each loop edge $f$ with $t(f)=o(f)=u$ replace $f_{\alpha}$
  with $gf_{\alpha}$ and $f_{\omega}$ with $f_{\omega}g^{-1}$.
\end{enumerate}
In this case we will say that \emph{$\mathcal B'$ is obtained from
  $\mathcal B$ by a folding move of type $A0$.}

If $u'\in B$, $u'\ne u$ is another vertex (whose vertex group is
therefore not changed by the move), we will say that this $A0$-move is
\emph{admissible with respect to $u'$.}
\end{defn}

\begin{figure}[here]
  \centerline{ \footnotesize \setlength{\unitlength}{.8cm}
\begin{picture}(16,4)
\put(12.5,2){\circle*{.07}}
\bezier{200}(12.5,2)(11.4,3)(10.6,3)
\bezier{200}(12.5,2)(11.4,1)(10.6,1)
\bezier{200}(10.6,3)(9.7,3)(9.7,2)
\bezier{200}(10.6,1)(9.7,1)(9.7,2)
\put(12.5,2){\line(4,3){2}}
\put(12.5,2){\line(4,-3){2}}
\put(14.5,3.5){\circle*{.07}}
\put(14.5,.5){\circle*{.07}}
\put(3.5,2){\circle*{.07}}
\put(6.5,2){\vector(1,0){2}}
\bezier{200}(3.5,2)(2.4,3)(1.6,3)
\bezier{200}(3.5,2)(2.4,1)(1.6,1)
\bezier{200}(1.6,3)(.7,3)(.7,2)
\bezier{200}(1.6,1)(.7,1)(.7,2)
\put(3.5,2){\line(4,3){2}}
\put(3.5,2){\line(4,-3){2}}
\put(5.5,3.5){\circle*{.07}}
\put(5.5,.5){\circle*{.07}}
\put(4.5,2.75){\vector(4,3){.01}}
\put(4.5,1.25){\vector(-4,3){.01}}
\put(13.5,2.75){\vector(4,3){.01}}
\put(13.5,1.25){\vector(-4,3){.01}}
\put(.7,2){\vector(0,1){.01}}
\put(9.7,2){\vector(0,1){.01}}
\put(5.1,3.6){$(B_3,v_3)$}
\put(14.1,3.6){$(B_3,v_3)$}
\put(5,.1){$(B_2,v_2)$}
\put(14,.1){$(B_2,v_2)$}
\put(3.8,1.88){$(B_1,v_1)$}
\put(12.8,1.88){$(gB_1g^{-1},v_1)$}
\put(3.2,3){$(a',e',b')$}
\put(12.05,3.05){$(ga',e',b')$}
\put(3.3,.7){$(a,e,b)$}
\put(12.1,.65){$(a,e,bg^{-1})$}
\put(-.8,2.05){$(\overline{a},\overline{e},\overline{b})$}
\put(7.6,2.25){$(g\overline{a},\overline{e},\overline{b}g^{-1})$}
\end{picture}}
\caption{A move of type $A0$ with $g\in A_{v_1}$}\label{fa0}
\end{figure}
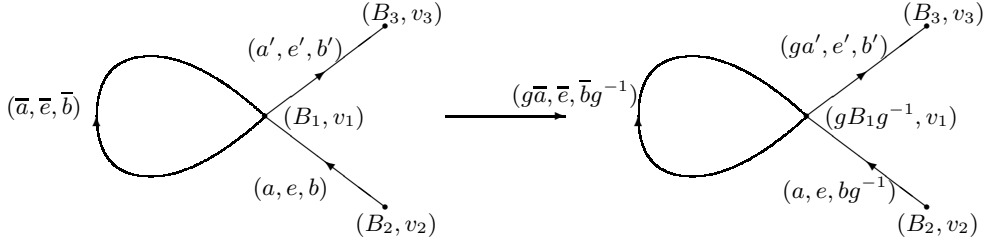

\begin{defn}[Bass-Serre move $A1$]
  Let $\mathcal B$ be an $\mathbb A$-graph. Suppose $f$ is an edge of
  $\mathcal B$ and that $c\in A_{[f]}$.

  Let $\mathcal B'$ be the $\mathbb A$-graph obtained from $\mathcal
  B$ by replacing $f_{\alpha}$ with $f_{\alpha}\alpha_{[e]}(c)^{-1}$
  and $f_{\omega}$ with $\omega_{[e]}(c)f_{\omega}$.

  In this case we will say that \emph{$\mathcal B'$ is obtained from
    $\mathcal B$ by a folding move of type $A1$.}
\end{defn}

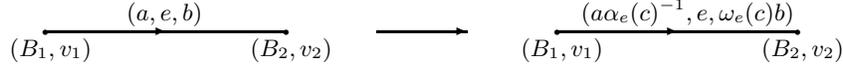
\begin{figure}[here]
  \centerline{ \footnotesize \setlength{\unitlength}{.8cm}
\begin{picture}(11,1)
\put(4.5,.5){\vector(1,0){1.5}}
\put(7.5,0.5){\line(3,0){4}}
\put(-1,0.5){\line(3,0){4}}
\put(-1,.5){\circle*{.07}}
\put(3,.5){\circle*{.07}}
\put(7.5,.5){\circle*{.07}}
\put(11.5,.5){\circle*{.07}}
\put(-1.6,.1){$(B_1,v_1)$}
\put(6.9,.1){$(B_1,v_1)$}
\put(2.4,.1){$(B_2,v_2)$}
\put(10.9,.1){$(B_2,v_2)$}
\put(.35,.7){$(a,e,b)$}
\put(7.9,.7){$(a\alpha_e(c)^{-1},e,\omega_e(c)b)$}
\put(9.5,.5){\vector(1,0){.01}}
\put(1,.5){\vector(1,0){.01}}
\end{picture}}
\caption{A move of type $A1$ with $c\in A_e$}\label{fa1}
\end{figure}

\begin{defn}[Simple adjustment $A2$]
  Let $\mathcal B$ be an $\mathbb A$-graph. Suppose $f$ is an edge of
  $\mathcal B$ and that $a'\in B_{o(f)}$.

  Let $\mathcal B'$ be the $\mathbb A$-graph obtained from $\mathcal
  B$ by replacing $f_{\alpha}$ with $a'f_{\alpha}$.

  In this case we will say that \emph{$\mathcal B'$ is obtained from
    $\mathcal B$ by a folding move of type $A2$.}
\end{defn}

\begin{figure}[here]
  \centerline{ \footnotesize \setlength{\unitlength}{.8cm}
\begin{picture}(11,1)
\put(4.5,.5){\vector(1,0){1.5}}
\put(7.5,0.5){\line(3,0){4}}
\put(-1,0.5){\line(3,0){4}}
\put(-1,.5){\circle*{.07}}
\put(3,.5){\circle*{.07}}
\put(7.5,.5){\circle*{.07}}
\put(11.5,.5){\circle*{.07}}
\put(-1.6,.1){$(B_1,v_1)$}
\put(6.9,.1){$(B_1,v_1)$}
\put(2.4,.1){$(B_2,v_2)$}
\put(10.9,.1){$(B_2,v_2)$}
\put(.35,.7){$(a,e,b)$}
\put(8.7,.7){$(a'a,e,b)$}
\put(9.5,.5){\vector(1,0){.01}}
\put(1,.5){\vector(1,0){.01}}
\end{picture}}
\caption{A move of type $A2$ with $a'\in B_1$}\label{fa2}
\end{figure}
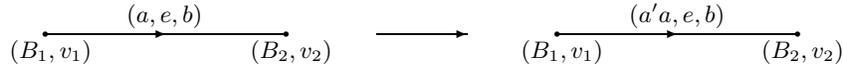

\subsection{Main Stallings type folding moves}\label{se:mainfold}

In this section we introduce folding moves that change the structure
of the underlying graph of an $\mathbb A$-graph. They can only be
applied to $\mathbb A$-graphs that are not folded. On the level of
underlying graphs these moves will correspond to the standard
Stallings folds.

\begin{conv}

  For the remainder of Section~\ref{se:mainfold} let $({\mathcal B},
  u_0)$ be an $\mathbb A$-graph with base vertex $u_0$.  Suppose
  $\mathcal B$ is not folded because case~(1) of
  Definition~\ref{defn:folded} applies. Thus there exist distinct
  edges $f_1$ and $f_2$ with $z=o(f_1)=o(f_2)$ and labels
  $(a_1,e,a_2)$ and $(a_2,e,b_2)$ such that $a_2=a'a_1\alpha_e(c)$ for
  some $c\in A_e$ and $a'\in B_z$.  Suppose further that $t(f_1)=x$
  and $t(f_2)=y$. Clearly $x$ and $y$ are of the same type $v\in VA$.
  We also denote the type of $z$ by $w\in VA$.

  \medskip By applying a move of type A2 to the edge $f_2$ we can
  change the label of $f_2$ to
  $({a'}^{-1}a_2,e,b_2)=({a'}^{-1}a'a_1\alpha_e(c),e,b_2)=(a_1\alpha_e(c),e,b_2)$.
  A move of type A1 then yields the label $(a_1,e,\omega_e(c)b_2)$ on
  $f_2$. We denote the resulting $\mathbb A$-graph by $\mathcal B'$.

\begin{figure}[here]
  \centerline{ \footnotesize \setlength{\unitlength}{.8cm}
\begin{picture}(16.5,7.2)
\put(1,5.5){\line(3,1){4}}
\put(1,5.5){\line(3,-1){4}}
\put(10,5.5){\line(3,1){4}}
\put(10,5.5){\line(3,-1){4}}
\put(6.5,5.5){\vector(1,0){2}}
\put(11.5,2.5){\vector(-2,-1){1.5}}
\put(7.3,5.6){A2}
\put(10.7,1.7){A1}
\put(1,5.5){\circle*{.07}}
\put(5,6.833){\circle*{.07}}
\put(5,4.166){\circle*{.07}}
\put(10,5.5){\circle*{.07}}
\put(14,6.833){\circle*{.07}}
\put(14,4.166){\circle*{.07}}
\put(0,5.1){$(B_z,w)$}
\put(9,5.1){$(B_z,w)$}
\put(5.1,6.8){$(B_x,v)$}
\put(14.1,6.8){$(B_x,v)$}
\put(5.1,4){$(B_y,v)$}
\put(14.1,4){$(B_y,v)$}
\put(2.2,6.55){$(a_1,e,b_1)$}
\put(11.2,6.55){$(a_1,e,b_1)$}
\put(1.4,4.1){$(a'a_1\alpha_e(c),e,b_2)$}
\put(10.5,4.1){$(a_1\alpha_e(c),e,b_2)$}
\put(3,6.166){\vector(3,1){.1}}
\put(12,6.166){\vector(3,1){.1}}
\put(3,4.833){\vector(3,-1){.1}}
\put(12,4.833){\vector(3,-1){.1}}
\put(4,1.5){\line(3,1){4}}
\put(4,1.5){\line(3,-1){4}}
\put(3,1.1){$(B_z,w)$}
\put(4,1.5){\circle*{.07}}
\put(8,2.833){\circle*{.07}}
\put(8,.166){\circle*{.07}}
\put(8.1,2.7){$(B_x,v)$}
\put(8.1,.05){$(B_y,v)$}
\put(5.2,2.55){$(a_1,e,b_1)$}
\put(4.6,.15){$(a_1,e,\omega_e(c)b_2)$}
\put(6,2.166){\vector(3,1){.1}}
\put(6,.833){\vector(3,-1){.1}}
\end{picture}}
\caption{Constructing $\mathcal B'$}\label{ff0}
\end{figure}
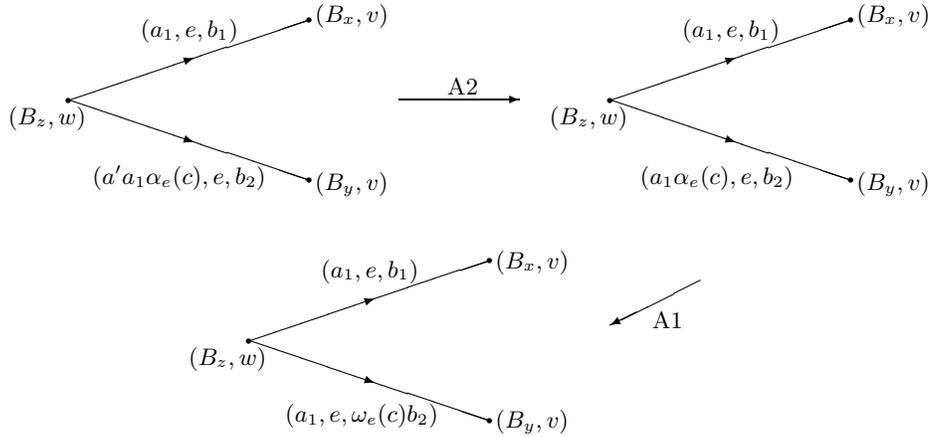

\end{conv}

\medskip

We will use $\mathcal B'$ as an intermediate object before defining
the main folding moves on $\mathcal B$. Note that $\mathcal B'$ is
obtained from $\mathcal B$ by moves that only alter labels of edges.
Thus $\mathcal B$ and $\mathcal B'$ have the same underlying graphs as
well as the same vertex groups.

It is possible that two or more of the vertices (that are drawn as
distinct vertices) coincide. To simplify notations we put $\bar
b_2:=\omega_e(c)b_2$. We now introduce four different types of folds,
$F1-F4$. They are distinguished by the topological type of the
subgraph $f_1\cup f_2$ in $B$. Each of these moves will be defined as
a sequence of several transformations, exactly one of which will
correspond to performing a Stallings fold identifying the edges $f_1$
and $f_2$ in $B$. That particular portion of a move $FN$,
$N=1,\dots,4$, will be called an \emph{elementary move of type}
$\overline{F}N$.

\begin{figure}[here]
  \centerline{ \footnotesize \setlength{\unitlength}{.8cm}
\begin{picture}(16.5,6)
\put(1,4.5){\line(3,1){4}}
\put(1,4.5){\line(3,-1){4}}
\put(10,4.5){\line(3,1){4}}
\put(10,4.5){\line(3,-1){4}}
\put(3.5,.5){\line(1,0){5}}
\put(6.5,4.5){\vector(1,0){2}}
\put(10,2.5){\vector(-2,-1){1.5}}
\put(7.3,4.6){A0}
\put(9.2,1.7){$\overline{F}1$}
\put(1,4.5){\circle*{.07}}
\put(5,5.833){\circle*{.07}}
\put(5,3.166){\circle*{.07}}
\put(10,4.5){\circle*{.07}}
\put(14,5.833){\circle*{.07}}
\put(14,3.166){\circle*{.07}}
\put(3.5,.5){\circle*{.07}}
\put(8.5,.5){\circle*{.07}}
\put(.2,4.15){$(B_z,w)$}
\put(9.2,4.15){$(B_z,w)$}
\put(2.3,.15){$(B_z,w)$}
\put(8.6,.15){$(\langle B_x, b_1^{-1}\bar b_2 B_y{\bar b}_2^{-1}b_1\rangle, v)$}
\put(5.1,5.8){$(B_x,v)$}
\put(14.1,5.8){$(B_x,v)$}
\put(5.1,3){$(B_y,v)$}
\put(13,2.75){$(b_1^{-1}\bar b_2 B_y{\bar b}_2^{-1}b_1,v)$}
\put(2.3,5.5){$(a_1,e,b_1)$}
\put(11.3,5.5){$(a_1,e,b_1)$}
\put(5.5,.8){$(a_1,e,b_1)$}
\put(2.3,3.3){$(a_1,e,\bar b_2)$}
\put(11.3,3.3){$(a_1,e,b_1)$}
\put(5.8,.5){\vector(1,0){.3}}
\put(3,5.166){\vector(3,1){.1}}
\put(12,5.166){\vector(3,1){.1}}
\put(3,3.833){\vector(3,-1){.1}}
\put(12,3.833){\vector(3,-1){.1}}
\end{picture}}
\caption{A move of type $F1$}\label{ff1}
\end{figure}
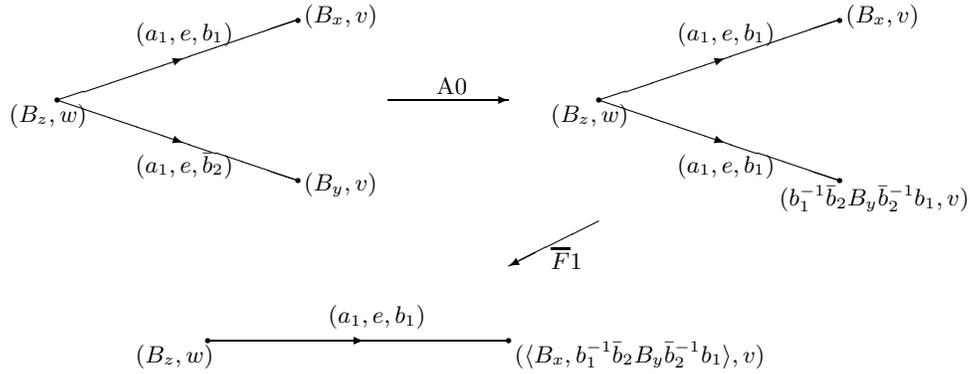

\begin{defn}[Simple fold $F1$]
  Suppose $f_1$ and $f_2$ are two distinct non-loop edges and that
  $t(f_1)\ne t(f_2)$. Possibly after exchanging $f_1$ and $f_2$ we can
  assume that $t(f_2)$ is not the base vertex $u_0$ of $\mathcal B$.

  We first perform a move of type $A0$ on $\mathcal B'$ at the vertex
  $t(f_2)=y$ making the label of $f_2$ to be $(a_1,e,b_1)$ and the
  label of $t(f_2)$ to be $(b_1^{-1}\bar b_2 B_y{\bar
    b}_2^{-1}b_1,v)$.  Now both $f_1$ and $f_2$ have label
  $(a_1,e,b_1)$.

  Next we identify the edges $f_1$ and $f_2$ into a single edge $f$
  with label $(a_1,e,b_1)$, as illustrated in Figure~\ref{ff1}.  The
  label of the vertex $t(f)$ is set to be
\[
(\langle B_x, b_1^{-1}\bar b_2 B_y{\bar b}_2^{-1}b_1\rangle, v).\] The
other labels do not change.

We call this last operation an {\em elementary move of type}
$\overline{F}1$ and say that the resulting $\mathbb A$-graph is
\emph{obtained from the original $\mathbb A$-graph $\mathcal B$ by a
  move of type $F1$}.
\end{defn}

\begin{figure}[here]
  \centerline{ \footnotesize \setlength{\unitlength}{.8cm}
\begin{picture}(17,6.5)
\put(5,3.666){\line(-1,0){3}}
\put(14,3.666){\line(-1,0){3}}
\put(6.5,5){\vector(1,0){2}}
\put(12,2.5){\vector(-3,-2){1.5}}
\put(3.5,3.666){\vector(1,0){.01}}
\put(2,6.2){\vector(1,0){.01}}
\put(12.5,3.666){\vector(1,0){.01}}
\put(11,6.2){\vector(1,0){.01}}
\put(8,2.7){\vector(1,0){.01}}
\put(7.3,5.1){A0}
\put(11.5,1.8){$\overline{F}2$}
\put(8,.166){\circle*{.07}}
\put(7.35,2.95){$(a_1,e,b_1)$}
\put(5,3.666){\circle*{.07}}
\put(2,3.666){\circle*{.07}}
\put(14,3.666){\circle*{.07}}
\put(11,3.666){\circle*{.07}}
\put(6,-0.25){$(\langle B_z, b_1^{-1}\bar b_2 B_y{\bar b}_2^{-1}b_1\rangle, v)$}
\put(4.8,3.3){$(B_y,v)$}
\put(1.3,3.3){$(B_z,v)$}
\put(10.3,3.3){$(B_z,v)$}
\put(13,3.25){$(b_1^{-1}\bar b_2 B_y{\bar b}_2^{-1}b_1,v)$}
\put(1.35,6.45){$(a_1,e,b_1)$}
\put(10.35,6.45){$(a_1,e,b_1)$}
\put(2.8,3.8){$(a_1,e,\bar b_2)$}
\put(11.8,3.8){$(a_1,e,b_1)$}
\bezier{200}(2,3.666)(1,4.3)(1,5.066)
\bezier{200}(2,3.666)(3,4.3)(3,5.066)
\bezier{200}(2,6.2)(1,6.2)(1,5.066)
\bezier{200}(2,6.2)(3,6.2)(3,5.066)
\bezier{200}(11,3.666)(10,4.3)(10,5.066)
\bezier{200}(11,3.666)(12,4.3)(12,5.066)
\bezier{200}(11,6.2)(10,6.2)(10,5.066)
\bezier{200}(11,6.2)(12,6.2)(12,5.066)
\bezier{200}(8,.166)(7,.8)(7,1.566)
\bezier{200}(8,.166)(9,.8)(9,1.566)
\bezier{200}(8,2.7)(7,2.7)(7,1.566)
\bezier{200}(8,2.7)(9,2.7)(9,1.566)
\end{picture}}
\caption{A move of type $F2$}\label{ff2}
\end{figure}
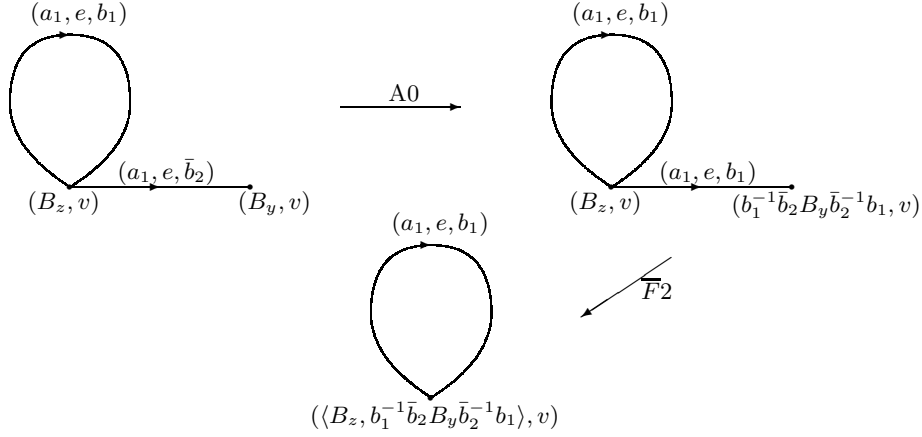

\begin{defn}[Mixed fold $F2$]
  Suppose now that $f_1$ is a loop edge and that $f_2\ne f_1$ is a
  non-loop edge. (The opposite situation is analogous).

  This implies that $e$ is a loop-edge in $A$ based at the vertex
  $v=w$.

  We first perform move $A0$ on $\mathcal B'$ making the label of
  $f_2$ to be $(a_1,e,b_1)$. Next we fold the edges $f_1$ and $f_2$
  into a single loop-edge $f$ with label $(a_1, e,b_1)$, as shown in
  Figure~\ref{ff2}. The label of $o(f)=t(f)$ is set to be
\[
(\langle B_z, b_1^{-1}\bar b_2 B_y{\bar b}_2^{-1}b_1\rangle, v).\] We
call this operation an {\em elementary move of type} $\bar F2$.

If $y=t(f_2)=u_0$ we then perform the auxiliary move $A0$
corresponding to the element ${\bar b_2}^{-1}b_1$.

We will say that the resulting $\mathbb A$-graph \emph{is obtained
  from $\mathcal B$ by a folding move of type $F2$.}
\end{defn}

\begin{figure}[here]
  \centerline{ \footnotesize \setlength{\unitlength}{.8cm}
\begin{picture}(15,3)
\put(3.5,1.5){\circle*{.07}}
\put(11.5,1.5){\circle*{.07}}
\put(7.9,1.5){\vector(1,0){2}}
\bezier{200}(3.5,1.5)(4.6,2.5)(5.4,2.5)
\bezier{200}(3.5,1.5)(4.6,.5)(5.4,.5)
\bezier{200}(5.4,2.5)(6.3,2.5)(6.3,1.5)
\bezier{200}(5.4,.5)(6.3,.5)(6.3,1.5)
\bezier{200}(3.5,1.5)(2.4,2.5)(1.6,2.5)
\bezier{200}(3.5,1.5)(2.4,.5)(1.6,.5)
\bezier{200}(1.6,2.5)(.7,2.5)(.7,1.5)
\bezier{200}(1.6,.5)(.7,.5)(.7,1.5)
\bezier{200}(11.5,1.5)(12.6,2.5)(13.4,2.5)
\bezier{200}(11.5,1.5)(12.6,.5)(13.4,.5)
\bezier{200}(13.4,2.5)(14.3,2.5)(14.3,1.5)
\bezier{200}(13.4,.5)(14.3,.5)(14.3,1.5)
\put(14.3,1.5){\vector(0,-1){.01}}
\put(6.3,1.5){\vector(0,-1){.01}}
\put(.7,1.5){\vector(0,-1){.01}}
\put(-1.1,1.2){$(a_1,e,b_1)$}
\put(14.4,1.2){$(a_1,e,b_1)$}
\put(6.35,1.1){$(a_1,e,\bar b_2)$}
\put(2.93,2){$(B_z,u)$}
\put(9.5,2.1){$(\langle B_z, b_1^{-1}\bar b_2\rangle,u)$}
\put(8.6,1){$\overline{F}3$}
\end{picture}}
\caption{A move of type $F3$}\label{ff3}
\end{figure}
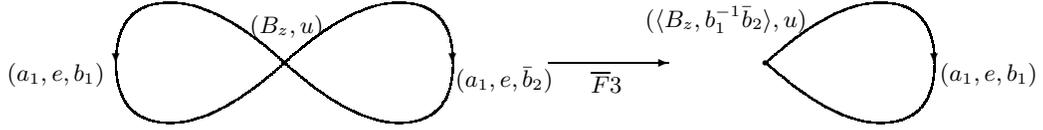

\begin{defn}[Loop fold $F3$]
  Suppose $f_1$ and $f_2$ are distinct loop-edges, so that $x=y=z$,
  $v=w$ and $e$ is a loop-edge at $v=w$ in $A$.

  We identify the edges $f_1$ and $f_2$ in $\mathcal B'$ into a single
  loop with label $(a_1, e,b_1)$, as shown in Figure~\ref{ff3}. The
  new label of $z$ is set to be \[ (\langle B_z, b_1^{-1}\bar
  b_2\rangle, v).\]

  We call this last operation an {\em elementary move of type}
  $\overline{F}3$ and say that the resulting $\mathbb A$-graph
  \emph{is obtained from $\mathcal B$ by a folding move of type $F3$.}
\end{defn}

\begin{defn}[Double-edge fold $F4$]
  Suppose that $f_1$ are both non-loop edges such that
  $x=t(f_1)=t(f_2)=y$.

  We identify the edges edges $f_1$ and $f_2$ of $\mathcal B'$ into a
  single edge $f$ with label $(a_1,e,b_1)$. We set the label of $t(f)$
  to be
\[
(\langle B_x, b_1^{-1}\bar b_2\rangle, v).\]

We call this last operation an {\em elementary move of type}
$\overline{F}4$ and say that the resulting $\mathbb A$-graph is
obtained from $\mathcal B$ by a {\em folding move of type $F4$.}

\begin{figure}[here]
  \centerline{ \footnotesize \setlength{\unitlength}{.8cm}
\begin{picture}(10,3)
\bezier{200}(-2,1.5)(-.33,2.5)(.5,2.5)
\bezier{200}(-2,1.5)(-.33,.5)(.5,0.5)
\bezier{200}(.5,2.5)(1.33,2.5)(3,1.5)
\bezier{200}(.5,0.5)(1.33,.5)(3,1.5)
\put(4.5,1.5){\vector(1,0){1.5}}
\put(.5,2.5){\vector(1,0){.01}}
\put(.5,.5){\vector(1,0){.01}}
\put(9.5,1.5){\vector(1,0){.01}}
\put(7.5,1.5){\line(3,0){4}}
\put(-2,1.5){\circle*{.07}}
\put(3,1.5){\circle*{.07}}
\put(7.5,1.5){\circle*{.07}}
\put(11.5,1.5){\circle*{.07}}
\put(5,1.6){$\overline{F}4$}
\put(-.3,2.65){$(a_1,e,b_1)$}
\put(-.3,0){$(a_1,e,\bar b_2)$}
\put(6.5,1.05){$(B_z,w)$}
\put(10.3,1.1){$(\langle B_x,b_1^{-1}\bar b_2\rangle,v)$}
\put(8.7,1.75){$(a_1,e,b_1)$}
\put(-2.8,1){$(B_z,w)$}
\put(3.1,1.15){$(B_x,v)$}
\end{picture}}
\caption{A move of type $F4$}\label{ff4}
\end{figure}
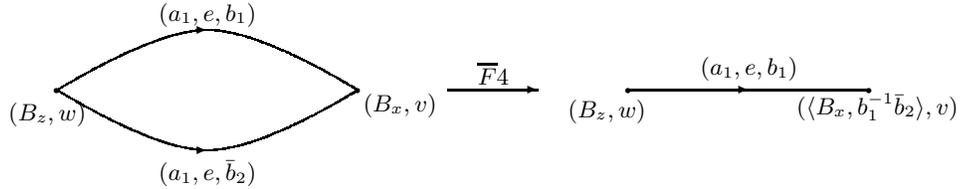

\end{defn}

\subsection{Edge-equalizing moves}

We will now introduce two folding moves that can be applied to an
$\mathbb A$-graph that is not folded because of the second condition
in Definition~\ref{defn:folded}.

Thus suppose $\mathcal B$ is an $\mathbb A$-graph with a base-vertex
$u_0$ and that there is an edge $f\in EB$ with label $(a,e,b)$, with
$o(f)$ labeled $(A',u)$ and $t(f)$ labeled $(B',v)$ such that
$\alpha_e^{-1}(a^{-1}A'a)\ne \omega_e^{-1}(bB'b^{-1})$.

\begin{defn}[Equalizing an edge group $F5$]
  Let $\mathcal B$ be an $\mathbb A$-graph.  Suppose $f$ is a non-loop
  edge of $\mathcal B$ with the label $(a, e, b)$. Let $(A',w)$ be the
  label of $z=o(f)$ and let $(A'',v)$ be the label of $t=t(f)$.

  Put $C:=\langle \alpha_e^{-1}
  (a^{-1}A'a),\omega_e^{-1}(bA''b^{-1})\rangle \le A_e$.

  Let $\mathcal B'$ be the $\mathbb A$-graph obtained from $\mathcal
  B$ by replacing the label of $z$ with the label $(\langle A', a
  \alpha_e(C) a^{-1}\rangle, w)$ and the label of $t$ with $(\langle
  A'', b^{-1} \omega_e(C) b\rangle,v)$.  In this case we will say that
  \emph{$\mathcal B'$ is obtained from $\mathcal B$ by a move of type
    $F5$}.
\end{defn}

\begin{figure}[here]
  \centerline{ \footnotesize \setlength{\unitlength}{.8cm}
\begin{picture}(11,1)
\put(4.5,.5){\vector(1,0){1.5}}
\put(5,.65){$F5$}
\put(-1,0.5){\line(3,0){4}}
\put(7.5,0.5){\line(3,0){4}}
\put(-1,.5){\circle*{.07}}
\put(3,.5){\circle*{.07}}
\put(7.5,.5){\circle*{.07}}
\put(11.5,.5){\circle*{.07}}
\put(-1.6,.1){$(A',w)$}
\put(6,.1){$(\langle A', a \alpha_e(C) a^{-1}\rangle,w)$}
\put(2.4,.1){$(A'',v)$}
\put(10.5,.1){$(\langle A'', b^{-1} \omega_e(C) b\rangle,v)$}
\put(.35,.7){$(a,e,b)$}
\put(8.8,.7){$(a,e,b)$}
\put(9.5,.5){\vector(1,0){.01}}
\put(1,.5){\vector(1,0){.01}}
\end{picture}}
\caption{A move of type $F6$ with $C:=\langle \alpha_e^{-1}
  (a^{-1}A'a),\omega_e^{-1}(bA''b^{-1})\rangle$}\label{ff5}
\end{figure}

\begin{defn}[Equalizing a loop-edge group $F6$]
  Suppose $f$ is a loop edge of $\mathcal B$ with the label $(a, e,
  b)$. Let $(A',v)$ be the label of $z=o(f)=t(f)$.

  Put $C:=\langle \alpha_e^{-1}
  (a^{-1}A'a),\omega_e^{-1}(bA'b^{-1})\rangle \le A_e$.

  Let $\mathcal B'$ be the $\mathbb A$-graph obtained from $\mathcal
  B$ by replacing the label of $z$ with
\[
(\langle A', a \alpha_e(C) a^{-1}, b^{-1} \omega_e(C) b\rangle, v).
\]
In this case we will say that \emph{$\mathcal B'$ is obtained from
  $\mathcal B$ by a move of type $F6$}.
\end{defn}

\begin{figure}[here]
  \centerline{ \footnotesize \setlength{\unitlength}{.8cm}
\begin{picture}(20,3)
\put(3.5,1.5){\circle*{.07}}
\put(11.5,1.5){\circle*{.07}}
\put(7.9,1.5){\vector(1,0){2}}
\bezier{200}(3.5,1.5)(4.6,2.5)(5.4,2.5)
\bezier{200}(3.5,1.5)(4.6,.5)(5.4,.5)
\bezier{200}(5.4,2.5)(6.3,2.5)(6.3,1.5)
\bezier{200}(5.4,.5)(6.3,.5)(6.3,1.5)
\bezier{200}(11.5,1.5)(12.6,2.5)(13.4,2.5)
\bezier{200}(11.5,1.5)(12.6,.5)(13.4,.5)
\bezier{200}(13.4,2.5)(14.3,2.5)(14.3,1.5)
\bezier{200}(13.4,.5)(14.3,.5)(14.3,1.5)
\put(14.3,1.5){\vector(0,-1){.01}}
\put(6.3,1.5){\vector(0,-1){.01}}
\put(14.4,1.2){$(a,e,b)$}
\put(6.4,1.2){$(a,e,b)$}
\put(2.4,1.65){$(A',v)$}
\put(10.3,1.65){$(A'',v)$}
\put(8.6,1){$F6$}
\put(7,0){$A''=\langle A', a \alpha_e(C) a^{-1}, b^{-1} \omega_e(C) b\rangle$}
\end{picture}}
\caption{A move of type $F6$ with $C:=\langle \alpha_e^{-1}
  (a^{-1}A'a),\omega_e^{-1}(bA'b^{-1})\rangle$}\label{ff6}
\end{figure}

Notice that each of the folding moves corresponds to a graph-morphism
between the underlying graphs which preserves types of vertices and
edges. In case of moves $F1-F4$ this morphism reduces the number of
edge-pairs by one. For moves $A0-A3,F5-F6$ the morphism is the
identity map. Moreover, the moves $F3,F4$ decrease the rank of the
fundamental group of the underlying graph $B$ by one, while $F1$ and
$F2$ do not change it.

The following important proposition states that folding moves preserve
the subgroup defined by an $\mathbb A$-graph.

\begin{Proposition}\label{fold_effect1}
  Let $\mathbb A$ be a graph of groups with a base-vertex $v_0$.
  Denote $G=\pi_1({\mathbb A}, v_0)$ and $X=\widetilde{({\mathbb A},
    v_0)}$.  Let $\mathcal B'$ be an $\mathbb A$-graph obtained from
  $\mathcal B$ by one of the folding moves $A0-A2, F1-F6$ (where an
  $A0$-fold is $u_0$-admissible).  Let $u_0$ be a vertex of $B$ and
  let $u_0'$ be the image of $u_0$ in $B'$. Suppose the type of the
  vertices $u_0, u_0'$ is $v_0\in VA$.

  Then there exists a canonical epimorphism $\gamma: \pi_1({\mathbb
    B}, u_0)\to \pi_1({\mathbb B'}, u_0')$ and a $\gamma$-equivariant
  simplicial map $\xi: \widetilde{({\mathbb B}, u_0)}\to
  \widetilde{({\mathbb B'}, u_0')}$ preserving the base-points such
  that the diagrams 
   $$\begin{array}{cccc}\xymatrix{
\pi_1({\mathbb B}, u_0)\ar[r]^{\gamma} \ar[dr]^{\nu} &\pi_1({\mathbb B'}, u_0')\ar[d]^{\nu'}
\\
&G}&\ \ {\begin{array}{c} \\ \\ \\\hbox{and}\end{array}}\ \ &\xymatrix{
\widetilde{({\mathbb B}, u_0)}\ar[r]^{\xi} \ar[dr]^{\phi} &\widetilde{({\mathbb B'}, u_0')}\ar[d]^{\phi'}
\\
&X}&\ \ {\begin{array}{c} \\ \\ \\\hbox{ commute.}\end{array}}\end{array}$$

%

Hence the images of $\phi: \widetilde{({\mathbb B}, u_0)}\to
\widetilde{({\mathbb A}, v_0)}$ and $\phi': \widetilde{({\mathbb B'},
  u_0')}\to \widetilde{({\mathbb A}, v_0)}$ coincide, and
\[
\overline {L({\mathcal B}, u_0)}=\overline {L({\mathcal B'}, u_0')}\le
G=\pi_1({\mathbb A}, v_0).
\]

\end{Proposition}

\begin{proof}
  The proof of this proposition is a straightforward exercise in
  Bass-Serre theory. We will sketch a sample argument for the folding
  move $F1$ and leave the other cases to the reader.

  If we assume that auxiliary moves have already been shown to satisfy
  Proposition~\ref{fold_effect1}, we can assume that the move $F1$ is
  actually an elementary move of type $\overline{F}1$, i.e. that both
  $f_1$ and $f_2$ have labels $(a,e,b)$ and their terminal vertices
  $t(f_1),t(f_2)$ have labels $(B_x,v)$ and $(B_y,v)$ accordingly.
  The folding move $\bar F1$ identifies $f_1$ and $f_2$ into a single
  edge $f$ with label $(a,e,b)$ and with the label of $t(f)$ equal
  $(\langle B_x\cup B_y\rangle , v)$

  Denote the folding graph-map in this move by $P:B\to B'$ so that
  $P(u_0)=u_0'$ and $P(f_1)=P(f_2)=f$. Note that by definition of
  $\overline{F}1$ for any edge $f'\in EB$ with $f\not\in \{f_1, f_2,
  f_1^{-1}, f_2^{-1}\}$ we have $P(f')=f'$. Also by construction for
  every vertex $u_1\in VB$ we have $B_{u_1}\le B'_{P(u_1)}$.  Thus the
  map $P$ gives rise to the obvious map $\gamma$ which takes a
  $\mathbb B$-path from $u_1\in VB$ to $u_2\in VB$ to a $\mathbb
  B'$-path from $P(u_1)$ to $P(u_2)$. It is easy to see that $\gamma$
  respects the $\sim$-equivalence relation and therefore factors
  through to a group homomorphism, also denoted by $\gamma$
\[
\gamma:\pi_1({\mathbb B}, u_0) \longrightarrow \pi_1({\mathbb B'},
u_0').
\]

The only nontrivial statement about the properties of $\gamma$ is to
check that $\gamma$ is in fact ``onto''. It suffices to show that a
generating set for $\pi_1({\mathbb B'}, u_0')$ provided by
Lemma~\ref{gen} lies in the image of $\gamma$. Since the edges
$f_1,f_2$ being folded by a move of type $\overline{F}1$ are non-loops
in $B$, we can choose a spanning tree $T$ in $B$ which contains both
$f_1$ and $f_2$. Then the graph $P(T)$ obtained from $T$ by
identifying $f_1$ and $f_2$ is clearly a spanning tree for $B'$.
Suppose $f\in EB'-P(T)$. Then $f$ is in fact an edge of $B$ which lies
outside of $T$. Hence
 \[
 \gamma([u_0,o(f)]_T f [t(f), u_0]_T)=[u_0',o(f)]_{P(T)} f [t(f),
 u_0']_{P(T)}=:s_f
\]
and so the generator $s_f$ of $\pi_1({\mathbb B'}, u_0')$ belongs to
the image of $\gamma$.

Assume now that $u'\in VB'$ is a vertex of $B'$. We need to show that
the set $[u_0',u']_{P(T)} B'_{u'} [u', u_0']_{P(T)}$ is contained in
the image of $\gamma$. If $u'\ne t(f)$ then by construction
$u'=P(u')\in VB$ is a vertex of $B$ with $B_{u'}=B'_{u'}$. In this
case $P([u_0,u']_T)$ is the $P(T)$-geodesic from $u_0'=P(u_0)$ to $u'$
and so
\[
[u_0',u']_{P(T)} B'_{u'} [u', u_0']_{P(T)}=\gamma([u_0,u']_TB_{u'}[u',
u_0']_T)\subseteq im (\gamma),
\]
as required.  Suppose next that $u'=t(f)=P(t(f_1))=P(t(f_2))$. We will
assume that $f$ is contained in the $P(T)$-geodesic from $u_0'$ to
$u'$ as the other case is similar.

Recall that by construction we have
\[
B'_{u'}=\langle B_1,B_2\rangle=\langle
B_{t(f_1)}, B_{t(f_2)}\rangle.
\]

Thus it suffices to show that for $i=1,2$ the set
\[
[u_0',u']_{P(T)} B_i [u', u_0']_{P(T)}\] is contained in the image of
$\gamma$.  Since $[u_0', u']_{P(T)}=P([u_0, t(f_1)]_T)=P([u_0,
t(f_2)]_T)$, it follows that
\[
[u_0',u']_{P(T)} B_i [u', u_0']_{P(T)}=\gamma([u_0,t(f_i)]_T B_i
[t(f_i), u_0]_T)\subseteq im (\gamma),
\]
as required. Thus indeed $\gamma: \pi_1({\mathbb B}, u_0)\to
\pi_1({\mathbb B'}, u_0')$ is surjective.

We will now define $\xi: \widetilde{({\mathbb B}, u_0)}
\longrightarrow \widetilde{({\mathbb B'}, u_0')}$. Suppose $x$ is a
vertex of $\widetilde{({\mathbb B}, u_0)}$. Thus $x$ has the form
$x=\overline{p} B_{u_1}$ for some vertex $u_1\in VB$ and some $\mathbb
B$-path $p$ from $u_0$ to $u_1$. Then put
$\xi(x):=\overline{\gamma(p)} B_{P(u_1)}$. It is not hard to see that
$\xi$ is well-defined on the vertex set of $\widetilde{({\mathbb B},
  u_0)}$ and that it preserves the adjacency relation for vertices.
Thus indeed we have constructed a simplicial map $\xi: \widetilde{({\mathbb B}, u_0)} \longrightarrow
\widetilde{({\mathbb B'}, u_0')}$, as promised. The equivariant
properties of $\xi$ easily follow from the description of $\xi$ and
$\gamma$ given above and from the explicit construction of the maps
$\nu$ and $\phi$ given earlier in the proof of
Proposition~\ref{basic}. We leave the details to the reader.
\end{proof}

\begin{Lemma} Let $\mathcal B$ be an $\mathbb A$-graph.
\begin{enumerate}
\item The $\mathbb A$-graph $\mathcal B$ is folded if and only if none
  of the moves $F1-F6$ apply. Moreover, if $\mathcal B$ is folded then
  any application of moves of type $A0-A2$ produces another folded
  graph.
\item Suppose that $\mathcal B$ is not folded and case (1) of
  Definition~\ref{defn:folded} occurs. Then a move of type $F1-F4$
  can be applied to $\mathcal B$.
\item Suppose that $\mathcal B$ is not folded and case (2) of
  Definition~\ref{defn:folded} occurs. Then a move of type $F5$ or $F6$ can be
  applied to $\mathcal B$.
\end{enumerate}
\end{Lemma}

\begin{proof}
  The statement of the lemma follows immediately from the definition
  of a folded graph and the section that
  introduces the folding moves.
\end{proof}

\section{Finding the induced splitting algorithmically}

In this section we describe an explicit procedure for finding an
induced splitting for a subgroup and give a set of sufficient
conditions which allow one to do this algorithmically.

\smallskip The following notion allows us to easily construct a
(usually non-folded) $\mathbb A$-graph for a subgroup $U$ of
$G=\pi_1({\mathbb A}, v_0)$ given by a generating set $S\subset G$ of
$U$.

\begin{defn}[Wedge]\label{defn:wedge}

  Let $\mathbb A$ be a graph of groups with a base-vertex $v_0$ and
  let $S\subset G=\pi_1({\mathbb A}, v_0)$. For each $s\in S$ we
  choose a reduced $\mathbb A$-path $p_s$ from $v_0$ to $v_0$ such
  that $\overline {p_s}=s$. Put $P_S=\{p_s|s\in S\}$ .

  We construct an $\mathbb A$-graph $\mathcal B_0$ as follows.  The
  underlying graph $B_0$ has base-vertex called $u_0$ of type $v_0$.
  For each path $p_s\in P_S$ of length at least $1$ we write $p_s$ as
  $p_s=a_0, e_1, a_1,\dots, e_k, a_k$ and attach at the vertex $u_0$ a
  circle subdivided into $k$ edges. We give the first $k-1$ of these
  edges labels $(a_0, e_1, 1), \dots , (a_{k-2},e_{k-1}, 1)$
  accordingly. We label the last edge of the circle by $(a_{k-1},e_k,
  a_k)$. This describes the underlying graph $B_0$ of $\mathcal B_0$
  with the obvious assignment of types for vertices and edges (Note
  that $B_0$ is either a single vertex or a wedge of circles).  Every
  vertex $u\in VB$ different from $u_0$ and of type $v\in VA$ is
  assigned the label $(1,v)$ (so that the corresponding vertex group
  is trivial).

  Note that for each $p_s\in P_S$ of length zero we have $s=\overline
  p_s\in A_{v_0}$.  We assign the vertex $u_0$ of $B_0$ label
  $(K,v_0)$, where
\[
K=\langle \{s\in S\ |\ |p_s|=0\}\rangle \le A_{v_0}.
\]
This completely describes the $\mathbb A$-graph $\mathcal B_0$. We
call such an $\mathbb A$-graph an $S$-\emph{wedge}.
\end{defn}

\begin{exmp}\label{exwedge}
  Suppose that $\mathbb A$ is the edge-of groups with edge pair
  $\{e,e^{-1}\}$ and $o(e)=v_0$ and $t(e)=v$ such that
  $A_{v_0}=F(a,b)$, $A_v=F(c,d)$, $A_e=\langle a^2=c^3\rangle$ and
  that the boundary monomorphisms are the inclusion maps. In
  particular we have $$G=\pi_1({\mathbb A},v_0)=F(a,b)\ast_{a^2=c^3}
  F(c,d).$$

  Suppose that $S=\{ s_1=a^4,s_2=b^2,
  s_3=c^3d^{10},s_4=d^{10}\}\subset G$. Clearly we have $s_1,s_2\in
  A_{v_0}$ and we can choose $p_{s_3}$ and $p_{s_4}$ as
  $p_{s_3}=1,e,c^3d^{10}, e^{-1},1$ and $p_{s_4}=1,e,d^{10}, e^{-1},1$,
  respectively. The diagram of the $S$-wedge then looks as follows:

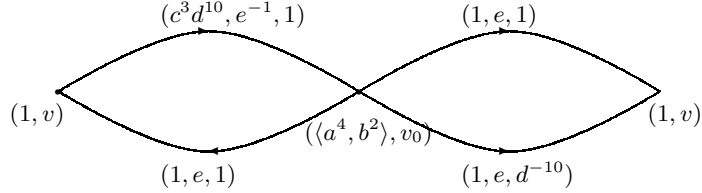
\begin{figure}[here]
  \centerline{ \footnotesize \setlength{\unitlength}{.8cm}
\begin{picture}(6,3)
\bezier{200}(-2,1.5)(-.33,2.5)(.5,2.5)
\bezier{200}(-2,1.5)(-.33,.5)(.5,0.5)
\bezier{200}(.5,2.5)(1.33,2.5)(3,1.5)
\bezier{200}(.5,0.5)(1.33,.5)(3,1.5)
\bezier{200}(3,1.5)(4.66,2.5)(5.5,2.5)
\bezier{200}(3,1.5)(4.66,.5)(5.5,0.5)
\bezier{200}(5.5,2.5)(6.33,2.5)(8,1.5)
\bezier{200}(5.5,0.5)(6.33,.5)(8,1.5)
\put(.5,2.5){\vector(1,0){.01}}
\put(.5,.5){\vector(-1,0){.01}}
\put(5.5,2.5){\vector(1,0){.01}}
\put(5.5,.5){\vector(1,0){.01}}
\put(-2,1.5){\circle*{.07}}
\put(3,1.5){\circle*{.07}}
\put(-.3,2.65){$(c^3d^{10}, e^{-1},1)$}
\put(-.3,0){$(1,e,1)$}
\put(4.7,2.65){$(1,e,1)$}
\put(4.7,0){$(1,e,d^{-10})$}
\put(-2.8,1){$(1,v)$}
\put(7.8,1){$(1,v)$}
\put(2.1,.7){$(\langle a^4,b^2\rangle,v_0)$}
\end{picture}}
\caption{The $S$-wedge of $S=\{ a^4,b^2, c^3d^{10}, d^{10}\}$}\label{s-wedge}
\end{figure}

\end{exmp}

\begin{Lemma}\label{wedge}
  Let $S,G,P_S, \mathbb B_0$ be as in the Definition~\ref{defn:wedge}
  and $U=\langle S\rangle$. Then
\[
\overline{L({\mathcal B}_0, u_0)}=U.\]

Moreover, the image of the map $$\phi: \widetilde{(B_0, u_0)} \to
X=\widetilde{({\mathbb A}, v_0)}$$
is equal to the tree $X_{U,x_0}$
(where $x_0$ is the base-vertex of $X$).

\end{Lemma}

\begin{proof}
  It is clear from the definitions that $S=\overline{P_S}\subseteq
  \overline{L({\mathcal B_0}, u_0)}=\nu(\pi_1({\mathbb B}, u_0))\le
  G$. Thus $U=\langle S\rangle\subset \overline{L({\mathcal B_0},
    u_0)}$. On the other hand Lemma~\ref{gen} implies that
  $\nu(\pi_1({\mathbb B}, u_0))$ is generated by $S$, and so
  $\overline{L({\mathcal B}_0, u_0)}=U$ as required.  Denote
  $H:=\pi_1({\mathbb B}, u_0)$.

  For each $s\in S$ with $|p_s|>0$ denote by $h_s$ the loop-path at
  $u_0$ in the wedge $B$ corresponding to $s$.  For each $s\in S$ with
  $|p_s|=0$ (so that $s\in A_{v_0}$) put $h_s=s$. Then each $h_s$ 
  defines a $\mathbb B$-path from $u_0$ to $u_0$ and hence an element
  of $H$ (if $|p_s|>0$ then $h_s$ gives a $\mathbb B$-path with
  trivial group elements inserted between the consecutive edges).
  Then $\nu(h_s)=s\in U\le G$ and $H=B_{u_0} \ast F(\{h_s \, | \,
  |p_s|>0\})$.  It follows from the definition of the Bass-Serre
  covering tree that the action of $H$ on $\widetilde{({\mathbb B},
    u_0)}$ is minimal, that is, it has no proper $H$-invariant
  subtrees. Denote the base-vertex of $\widetilde{({\mathbb B}, u_0)}$
  by $y_0$.

  Since $\pi_1({\mathbb B}, u_0)$ is generated by the set $\{h_s |
  s\in S\}$, it follows that
\[
\widetilde{({\mathbb B}, u_0)}=\cup_{h\in H} [y_0, h y_0]= H \big(
\cup_{s\in S} [y_0, h_s y_0]\big)
\]
and hence by Proposition~\ref{basic}
\[
\phi(\widetilde{(B_0, u_0)})=U \big( \cup_{s\in S} \phi([y_0, h_s
y_0])\big)
\]
Since $p_s$ is an $\mathbb A$-reduced path for each $s\in S$, we have
$\phi([y_0, h_s y_0])=[x_0, s x_0]$. Since $U=\langle S \rangle \le
G$, this implies
\[
U \big( \cup_{s\in S} \phi([y_0, h_s y_0])\big)=U\big( \cup_{s\in S}
[x_0,s x_0] \big)=\cup_{u\in U} [x_0,ux_0]=X_{U,x_0}
\]
as required.
\end{proof}

The following statement is an immediate corollary of
Proposition~\ref{prop:induced} and Proposition~\ref{fold_effect1}.

\begin{Proposition}[Abstract Folding Algorithm]\label{abstract}

  Let $\mathbb A$ be a graph of groups with a base-vertex $v_0$ and
  let $G=\pi_1({\mathbb A}, v_0)$. Suppose that $S\subset G$ and that
  $U=\langle S\rangle$.

  We first construct an $S$-wedge $\mathcal B_0$ with base-vertex
  $u_0$ as described in Definition~\ref{defn:wedge}. If this $\mathbb
  A$-graph is not folded, we start performing a sequence of folding
  moves $F1-F6$ (in an arbitrarily chosen order) and construct a
  sequence of based $\mathbb A$-graphs $({\mathcal B}_n, u_n)$, where
  each ${\mathcal B}_{n+1}$ is obtained from ${\mathcal B}_n$ by one of the folding moves
  $F1-F6$.

  If this sequence terminates in finitely many steps with a folded
  $\mathbb A$-graph $\mathcal B_m$, then we have $\overline{L({\mathcal B}_m,
    u_m)}=U$ and $\mathbb B_m$ gives the induced splitting for $U$ as
  described in Proposition~\ref{prop:induced}. \hfill$\Box$
\end{Proposition}

\begin{conv}
  When talking about actual algorithms related to graphs of groups we
  will not distinguish between an
  element of a vertex group and a word in the generators of that
  group. Thus, for example, when saying that we ``construct an
  $\mathbb A$-path $a_0, e_1, a_1, e_2 ,\dots$'' we actually mean
  constructing a sequence $w_0, e_1, w_1, e_2,\dots $ where $w_i$ is a
  word in the generators of the corresponding vertex group
  representing the element $a_i$. Moreover, we will assume that vertex
  and edge groups are explicitly given by recursive presentations on
  finite generating sets and that boundary monomorphisms are explicitly
  given by specifying the images of the generators of vertex groups in
  the appropriate edge groups.
\end{conv}

We will now describe a set of sufficient conditions which allows one
to algorithmically carry out the abstract procedure described in
Proposition~\ref{abstract}.

\begin{defn}\label{benign}

  We will say that a finite connected graph of finitely generated
  groups $\mathbb A$ is \emph{benign} if the following conditions are
  satisfied:

\begin{enumerate}
\item For each vertex $v\in VA$ and an edge $e\in EA$ with $o(e)=v$
  there is an algorithm with the following property. Given a finite
  set $X\subseteq A_v$ and an element $a\in A_v$ the algorithm decides
  whether $I=\langle X\rangle \cap a\alpha_e(A_e)$ is empty. If
  $I\neq\emptyset$, the algorithm produces an element of $I$.

\item Every edge group $A_e$ of $\mathbb A$ is {\em Noetherian}, that
  is, it contains no infinite ascending sequence of subgroups. (Being
  Noetherian is equivalent to saying that all subgroups are finitely
  generated).

\item Every edge group $A_e$ of $\mathbb A$ has solvable uniform
  membership problem, i.e. there is an algorithm which, given a finite
  subset $X\subseteq A_e$ and an element $a\in A_e$ decides whether or
  not $a\in \langle X\rangle$.
\item For each vertex $v\in VA$ and edge $e\in EA$ with $o(e)=v$ there
  is an algorithm with the following property. For any finite subset
  $X\subseteq A_v$ the algorithm computes a finite generating set for
  the subgroup $\alpha_e(A_e)\cap \langle X\rangle$.
\end{enumerate}

\end{defn}

\begin{rem}\label{normal}
  Notice that if $\mathbb A$ is benign, then $A_v$ has solvable
  membership problem with respect to $\alpha_e(A_e)$ (where $v=o(e),
  e\in EA$). Indeed, if $a\in A_v$ then $a\in \alpha_e(A_e)$ if and
  only if the intersection $\{1\} \cap a \alpha_e(A_e)$ is nonempty.
\end{rem}

\begin{Theorem}\label{construct}
  Let $\mathbb A$ be a benign graph of groups with base-vertex $v_0$.
  Then:

\begin{enumerate}
\item[(a)] There is an algorithm which, given a finite set in
  $G=\pi_1({\mathbb A}, v_0)$ generating a subgroup $U\le G$,
  constructs a folded $\mathbb A$-graph $\mathcal B$ with base-vertex
  $u_0$ such that $\overline{L({\mathcal B}, u_0)}=U$. In $\mathcal B$,
  each vertex group $B_u$ is given by its finite generating set of
  words in the generators of $A_{[v]}$.

  Moreover, $\nu:\pi_1({\mathbb B}, u_0)\to U\le G$ is an isomorphism,
  the map
  $$\phi: \widetilde{({\mathcal B}, u_0)}\to
  (X,x_0)=\widetilde{({\mathbb A} , v_0)}$$
  is injective and the image
  of $\phi$ is the tree $X_U=\cup_{u\in U} [x_0,ux_0]$.

\item[(b)] Suppose, in addition, that for each $v\in VA$ there is an
  algorithm which, given a finite subset $Y$ of $A_v$, produces a
  finite presentation for the subgroup of $A_v$ generated by $Y$ (thus each $A_v$ is coherent). Then there is an algorithm which,
  given a finite set $S\subseteq G$, constructs a finite presentation
  for the subgroup $U=\langle S\rangle \le G$.
\end{enumerate}
\end{Theorem}

\begin{rem}
  Thus by Proposition~\ref{prop:induced} the identification of $U$ with
  $\pi_1({\mathbb B}, u_0)$ via $\nu$ gives the induced splitting for
  $U\le G=\pi_1({\mathbb A} , v_0)$.  In particular, this
  identification gives us an explicit finite description of $U$,
  meaning that in $\mathcal B$ for each vertex $u\in VB$ of type $v\in
  VA$ the group-label of $u$ is given in the form $\langle X\rangle$,
  where $X$ is a finite subset of~$A_v$.
\end{rem}

\noindent{\em Proof of Theorem~\ref{construct}.} Let $S$ be a finite generating set of $U$. As $\mathbb A$ is benign we can find for any $s\in S$ a reduced $\mathbb A$-path $p_s$. Thus we can apply the
abstract folding algorithm as described in Proposition~\ref{abstract}.

We have to show that the process terminates in a finite number of
steps and that each step can be performed effectively.

Recall that by construction since $S$ is finite, the underlying graph
$B_0$ of $\mathcal B_0$ is finite. Moreover the vertex groups in
$\mathcal B_0$ are trivial with the possible exception of the
base-vertex $w_0\in VB_0$. By construction, the vertex group at $w_0$
is given by a finite generating set of cardinality at most $\#S$.

We can argue inductively that at each stage of the process for every
vertex $u\in VB_n$ of type $v\in VA$ the group $(B_n)_u\le A_v$ is
given by its finite generating set contained in $A_v$. At each stage
it is easy to decide whether $\mathcal B_n$ is folded.  Namely,
condition~(1) of Definition~\ref{benign} allows us to decide if Case
(1) of Definition~\ref{defn:folded} occurs.  Conditions (3) and (4) of
Definition~\ref{benign} allow us to decide if Case~(2) of
Definition~\ref{defn:folded} applies to $\mathcal B_n$. If $\mathcal
B_n$ turns out to be not folded, we perform one of the folded moves
$F1-F6$, whichever is appropriate.

By definition, performing folds of type $F1-F4$ allows us to
effectively represent the vertex groups of $\mathcal B_{n+1}$ by their
finite generating sets. Suppose now that $\mathcal B_{n+1}$ is
obtained from $B_{n}$ by a move of type $F5$ or $F6$.  Recall that
edge groups of $\mathbb A$ are Noetherian. Conditions (3) and (4) from
the definition of a benign graph of groups and the definitions of
folding moves $F5-F6$ allow us to effectively compute finite
generating sets for the vertex groups of $\mathcal B_{n+1}$.

Suppose that the sequence $({\mathcal B_n})$ is infinite. Each of the
moves of type $F1-F4$ reduces the number of edges in $\mathcal B_n$
and so can happen only finitely many types. Thus after a certain stage
only the moves of type $F5-F6$ (which do not change the underlying
finite graph) apply. Hence there is an edge to which moves of type
$F5-F6$ apply infinitely often. Each such move increases the edge
group of the corresponding edge in ${\mathbb B}_n$. This produces a
strictly increasing infinite sequence of subgroups in an edge-group of
$\mathbb A$, contradicting our assumption that edge-groups in $\mathbb
A$ are Noetherian.

Thus the sequence $\mathcal B_n$ terminates in finitely many steps
with a folded $\mathbb A$-graph $\mathcal B_m$, as required and part
(a) of Theorem~\ref{construct} is proved.

Once $\mathbb B$ as in Theorem~\ref{construct} is constructed, each
vertex (edge) group of $\mathbb B$ is given as a subgroup of some
vertex (edge) group of $\mathbb A$ generated by a given finite set of
elements. If the additional assumptions on $\mathbb A$ from part (b)
of Theorem~\ref{construct} hold, then we can recover finite
presentations for each vertex group of $\mathbb B$ and hence a finite
presentation for $U=\pi_1({\mathbb B}, w_0)$.~\hfill$\Box$

\begin{exmp}
  It is easy to produce an example of a non-benign graph of groups,
  where the folding algorithm described above does not necessarily
  terminate. For example, consider the HNN-extension of a free group
  $F=F(a,b)$ along the endomorphism $\phi:F\to F$, $\phi(a)=ab^2a,
  \phi(b)=ba^2b$:
  $$G=\langle a,b, e\, | \, e^{-1} ae=ab^2a, e^{-1}be=ba^2b\rangle=
  \langle a,b, e\, | \, e^{-1}\alpha_e(f)e=\omega_e(f), f\in F\rangle
  $$
  where $\alpha_e=Id_F$ and $\omega_e=\phi$. Thus we may think of
  $G$ as the fundamental group of the graph of groups $\mathbb A$
  consisting of a single vertex $v$, a single edge $e$ with
  $A_v=A_e=F$ and $\alpha_e=Id_F$ and $\omega_e=\phi$. The group $G$
  is torsion-free and word-hyperbolic~\cite{Ka00} by the Combination
  Theorem of Bestvina-Feighn~\cite{BF92}. Since $[a,e]\ne 1$ in $G$,
  there is $m>0$ such that $H=\langle e,a^m\rangle \le G$ is free of
  rank two. It is not hard to see that in this case $H\cap F$ is not
  finitely generated.  In fact $H\cap F$ is freely generated by the
  elements $e^{-i} a^m e^i=\phi^i(a^m)$. We can start the folding
  algorithm for $H$ with an $\mathbb A$-graph $\mathcal B$ consisting
  of a single vertex $u$ of type $v$, a single edge $f$ of type $e$
  with label $(1,e,1)$ and with $B_u=\langle a^m\rangle$. Then the
  folding algorithm results in a repeated application of an F6-move
  (no other moves are applicable) and produces an infinite sequence of
  $\mathbb A$-graphs $\mathcal B_0=\mathcal B, \mathcal B_1, \mathcal
  B_2, \dots$. The only difference between $\mathcal B_i$ and
  $\mathcal B$ is that in $\mathcal B_i$ the vertex group is $\langle
  a^m, \phi(a^m), \dots, \phi^i(a^m)\rangle$. This difficulty is
  caused by the fact that the edge-group in $\mathbb A$ is not
  Noetherian.

  A similar effect occurs in the direct product $F(a,b)\times \langle
  t\rangle$, thought of as an HNN-extension of $F(a,b)$ along the
  identity map, when we look at the subgroup $H=\langle tb, a\rangle$.

\begin{exmp}
  We illustrate the folding algorithm in Figure~\ref{ill}.  We start
  with the $S$-wedge discussed in Example~\ref{exwedge}. (Note that we
  have changed the orientation of two edges of the diagram.)

  The F5-move corresponds to the fact that $c^3=a^2$ and hence $a^2$
  has to be ``added" to the subgroup $\langle a^4,b^2\rangle$ yielding
  $\langle a^2,a^4,b^2\rangle=\langle a^2,b^2\rangle$. The final
  folded $\mathbb A$-graph corresponds to the induced splitting of
  $H=\langle S\rangle$ as $$H=\langle a^2,b^2\rangle*_{a^2=c^3}\langle
  c^3,d^{10}\rangle.$$
\end{exmp}
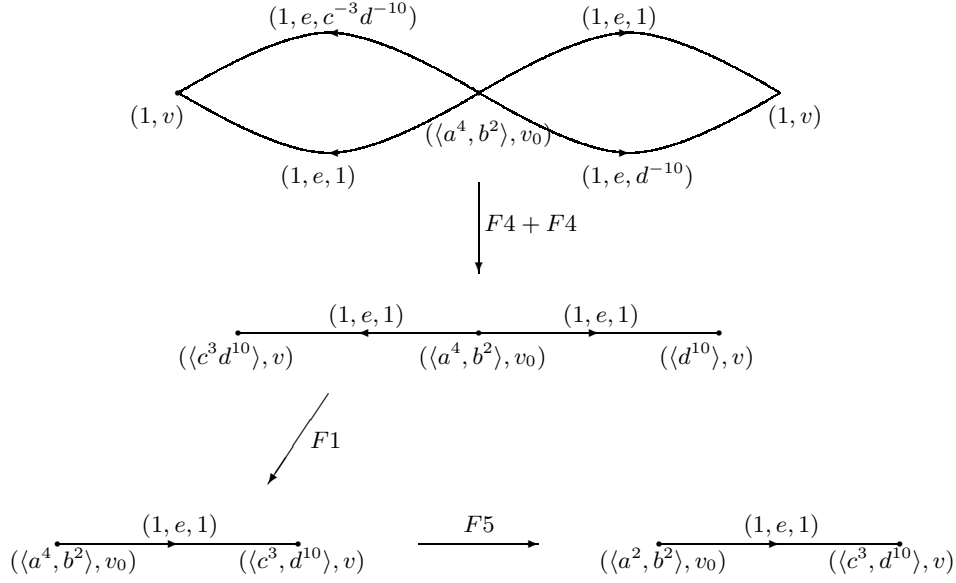
\begin{figure}[here]
  \centerline{ \footnotesize \setlength{\unitlength}{.8cm}
\begin{picture}(6,9.5)
\bezier{200}(-2,8)(-.33,9)(.5,9)
\bezier{200}(-2,8)(-.33,7)(.5,7)
\bezier{200}(.5,9)(1.33,9)(3,8)
\bezier{200}(.5,7)(1.33,7)(3,8)
\bezier{200}(3,8)(4.66,9)(5.5,9)
\bezier{200}(3,8)(4.66,7)(5.5,7)
\bezier{200}(5.5,9)(6.33,9)(8,8)
\bezier{200}(5.5,7)(6.33,7)(8,8)
\put(.5,9){\vector(-1,0){.01}}
\put(.5,7){\vector(-1,0){.01}}
\put(5.5,9){\vector(1,0){.01}}
\put(5.5,7){\vector(1,0){.01}}
\put(-2,8){\circle*{.07}}
\put(3,8){\circle*{.07}}
\put(-.5,9.15){$(1,e,c^{-3}d^{-10})$}
\put(-.3,6.5){$(1,e,1)$}
\put(4.7,9.15){$(1,e,1)$}
\put(4.7,6.5){$(1,e,d^{-10})$}
\put(-2.8,7.5){$(1,v)$}
\put(7.8,7.5){$(1,v)$}
\put(2.1,7.2){$(\langle a^4,b^2\rangle,v_0)$}
\put(3,6.5){\vector(0,-1){1.5}}
\put(3.1,5.7){$F4+F4$}
\put(-1,4){\line(1,0){8}}
\put(5,4){\vector(1,0){.01}}
\put(1,4){\vector(-1,0){.01}}
\put(-1,4){\circle*{.07}}
\put(3,4){\circle*{.07}}
\put(7,4){\circle*{.07}}
\put(.5,3){\vector(-2,-3){1}}
\put(.2,2.1){$F1$}
\put(-4,.5){\line(1,0){4}}
\put(6,.5){\line(1,0){4}}
\put(-2,.5){\vector(1,0){.01}}
\put(8,.5){\vector(1,0){.01}}
\put(-4,.5){\circle*{.07}}
\put(0,.5){\circle*{.07}}
\put(6,.5){\circle*{.07}}
\put(10,.5){\circle*{.07}}
\put(2,.5){\vector(1,0){2}}
\put(2.75,.7){$F5$}
\put(.5,4.2){$(1,e,1)$}
\put(4.4,4.2){$(1,e,1)$}
\put(-2.6,.7){$(1,e,1)$}
\put(7.4,.7){$(1,e,1)$}
\put(2,3.5){$(\langle a^4,b^2\rangle,v_0)$}
\put(-2,3.5){$(\langle c^3d^{10}\rangle,v)$}
\put(6,3.5){$(\langle d^{10}\rangle,v)$}
\put(-4.8,.1){$(\langle a^4,b^2\rangle,v_0)$}
\put(-1,.1){$(\langle c^3,d^{10}\rangle,v)$}
\put(5,.1){$(\langle a^2,b^2\rangle,v_0)$}
\put(8.8,.1){$(\langle c^3,d^{10}\rangle,v)$}
\end{picture}}
\caption{The folding algorithm applied to the $S$-wedge of Example~\ref{exwedge}}\label{ill}
\end{figure}
\end{exmp}

We will need the following simple observation which says that reduced
$\mathbb B$-paths of a folded $\mathbb A$-graph $\mathcal B$ are the
reduced $\mathbb A$-paths of the elements of $U=\overline{L({\mathcal
    B}, u_0)}$ up to the degree of freedom spelled out in the normal
form theorem for fundamental groups of graphs of groups. It provides a
criterion to decide whether an element lies in a subgroup represented
by a folded $\mathbb A$-graph.

\begin{Lemma}\label{belong} Let $\mathbb A$ be a graph of groups. Let $\mathcal
  B$ be a folded $\mathbb A$-graph. Suppose that $p=a_0, e_1,
  a_1,\dots, e_k, a_k$ is a reduced $\mathbb A$-path, where $k\ge 0$.
  Then $\overline{p}=\overline{\mu (q)}$ for some reduced $\mathbb
  B$-path $q$ if and only if there exists a reduced $\mathbb B$-path
\[
q=b_0,f_1,b_1, \dots, b_{k-1}, f_k, b_k\] and a sequence $c_i\in
A_{e_i}$, $i=1,\dots, k$ such that $[f_i]=e_i$ and

\begin{align*}
  &a_0=b_0 (f_1)_{\alpha}\alpha_{e_1}(c_1),
  \\
  &a_i=\omega_{e_i}(c_{i})^{-1}(f_i)_{\omega} b_i
  (f_{i+1})_{\alpha}\alpha_{e_{i+1}}(c_{i+1}) \text{ for } i=1,\dots,
  k-1\text{ and }\\
  & a_k=\omega_{e_k}(c_k)^{-1}(f_k)_{\omega} b_k
\end{align*}

\end{Lemma}

\begin{proof} The existence of a path $q$ and of
  $(b_i)_i$, $(c_i)_i$ with the required properties clearly implies
  that $\overline{p}=\overline{\mu(q)}$.

  If $\overline{p}=\overline{\mu(q)}$ for some reduced $\mathbb
  B$-path $q$, then the assertion follows from the normal form theorem
  applied to the product
  $p\mu(q)^{-1}$ which is trivial in $\pi_1({\mathbb A}, o(e_1))$.
\end{proof}

\begin{Theorem}\label{member}
  Let $\mathbb A$ be a benign graph of groups. Suppose also that each
  vertex group of $\mathbb A$ has solvable uniform membership problem.
  Let $v_0\in VA$ and denote $G=\pi_1({\mathbb A}, v_0)$.  Then the
  uniform membership problem for $G$ is solvable. That is to say,
  there is an algorithm which, given finitely many elements
  $h_1,\dots, h_k, g\in G$, decides whether $g$ belongs to the
  subgroup $H=\langle h_1,\dots, h_k\rangle\le G$.
\end{Theorem}

\begin{proof}

  Denote $S=\{ h_1,\dots, h_k\}$ and $U=\langle S\rangle\le G$. First
  we apply Theorem~\ref{construct} and construct a finite folded
  $\mathbb A$-graph $\mathcal B$ with a base-vertex $u_0$ such that
  $U=\overline{L({\mathcal B}, u_0)}$.  Every vertex group in
  $\mathcal B$ is given by its finite generating set contained in the
  appropriate vertex group of $\mathbb A$.

  Next we write $g$ as a reduced $\mathbb A$-path $p'$ from $v_0$ to
  $v_0$.  This is possible since by Remark~\ref{normal} every vertex
  group in $\mathbb A$ has solvable membership problem with respect to
  incident edge groups. Then $g\in U$ if and only if there exists a
  reduced $\mathbb B$-path $q'$ from $u_0$ to $u_0$ such that
  $\overline{p'}=\overline{\mu(q')}$. The assertion of the theorem now
  immediately follows from:

\medskip \noindent{\bf Claim.} There is an algorithm which, given a reduced $\mathbb A$-path
\begin{equation*}
p=a_0, e_1, a_1,\dots, e_k, a_k\tag{$*$}
\end{equation*}
from some vertex $v\in VA$ (possibly distinct from $v_0$) to $v_0$,
and given a vertex $u\in VB$,
decides if there exists a reduced $\mathbb B$-path
\[
q=b_0,f_1,b_1,\ldots ,b_{k-1},f_k,b_k
\]
from $u\in VB$ to $u_0$ such that
$\overline{\mu(q)}=\overline{p}$.

\medskip We prove the Claim by induction on the length $k$ of $p$. For $k=0$ the Claim is equivalent to deciding, given an
element $a_0\in A_{v_0}$, whether $a_0\in B_{u_0}$. This is possible
since $B_{u_0}\subset A_{v_0}$ is a finitely generated subgroup and by
assumption the group $A_{v_0}$ has solvable uniform membership
problem.

\smallskip Suppose now that $k>0$ and the algorithm exists for
reduced $\mathbb A$-paths of length $k-1$.

\smallskip Let $p$ be a reduced $\mathbb A$-path of length $k$ from $v$ to $v_0$
as in ($*$). If a path $q$ as in the Claim exists then it follows from Lemma~\ref{belong} that there is such
a path $q$ with the property that $[f_1]=e_1$ and
$a_0=b_0(f_1)_{\alpha}\alpha_{e_1}(c_1)$ for some $c_1\in A_{e_1}$.

\smallskip Observe first that we can decide whether there exists an
edge $f$ with $[f]=e_1$, $b\in B_{u}$ and $c\in A_{e_1}$ such that
$a_0=bf_{\alpha}\alpha_{e_1}(c)$ and can find them if they do exist. (If there are no such $f,b,c$ then
by the previous remark the required $q$ does not exist).

\smallskip Since there are only finitely many edges in $\mathcal B$ of type $e_1$
emanating at $u$, we may assume that we are dealing with a fixed edge
$f$ and looking for $b$ and $c$ with the above properties. Recall that
$a_0, f_{\alpha}\in A_{v}$ are given. Thus we want to know if there
are $b\in B_{u}$ and $c\in A_{e_1}$ such that $a_0=bf_{\alpha}\alpha_{e_1}(c)$, i.e. such that $f_{\alpha}^{-1}bf_{\alpha}=(f_{\alpha}^{-1}a_0)\alpha_{e_1}(c^{-1})$. Thus the existence of such $b$ and $c$ is equivalent to $f_{\alpha}^{-1}B_{u}f_{\alpha}\cap
(f_{\alpha}^{-1}a_0) \alpha_{e_1}(A_{e_1})\ne \emptyset$. This can be
checked by condition (1) of Definition~\ref{benign} since $\mathbb A$
is benign. Moreover, condition~(1) of Definition~\ref{benign} allows
us to find such (not necessarily unique) $b$ and $c$ if they exist.

\smallskip Suppose now that we have found $f,b$ and $c$ as above, so
that $a_0=bf_{\alpha}\alpha_{e_1}(c)$.

We now observe that if $q$ as in the Claim exists, then there is such
a $q$ with $b_0=b$. Indeed, if $q$ is as in the Claim then by
Lemma~\ref{belong} $[f_1]=e_1$ and there exists an element $c_1\in
\alpha_{e_1}(A_{e_1})$ such that
$a_0=b_0(f_1)_{\alpha}\alpha_{e_1}(c_1)=bf_{\alpha}\alpha_{e_1}(c)$.
The assumption that $\mathcal B$ is folded implies that $f=f_1$. Hence
$f_{\alpha}=(f_1)_{\alpha}$ and
\[
b^{-1}b_0=f_{\alpha}\alpha_{e_1}(cc_1^{-1})f_{\alpha}^{-1}\in
f_{\alpha}\alpha_{e_1}(A_{e_1})f_{\alpha}^{-1}.\] As $\mathcal B$ is folded this implies that $b^{-1}b_0\in \alpha_f(B_f)$. It follows that the
$\mathbb B$-path $q$ is equivalent to a $\mathbb B$-path starting with $b$,
as required. We denote this new $\mathbb B$-path again by $q$.

As $f_1=f$, $b_0=b$ in $q$ and $\overline{\mu(q)}=\overline{p}$ it follows that 
\[
\overline{b f_{\alpha} e_1 f_{\omega} b_1 (f_2)_{\alpha} e_2
  \dots}=\overline{a_0e_1a_1e_2\dots}=\overline{b f_{\alpha}
  \alpha_{e_1}(c) e_1 a_2 e_2\dots}=\overline{b
  f_{\alpha}e_1\omega_{e_1} (c) a_2e_2 \dots}
\]

and hence

\[
\overline{b_1 (f_2)_{\alpha}e_2\dots}=\overline{(f_{\omega})^{-1}
  \omega_{e_1} (c) a_1 e_2a_2\dots}.
\]

Thus to decide if a desired $q$ exists we need to determine if for the
reduced $\mathbb A$-path

\[
p'=(f_{\omega})^{-1} \omega_{e_1} (c) a_1, e_2, a_2, \ldots
,a_{k-1},e_k, a_k
\]
from $t(e_1)$ to $v_0$ in $\mathbb A$ there exists a path $q'$
starting at $t(f_1)$ as in
the Claim. This is possible by the inductive hypothesis since
$|p'|=k-1$.
\end{proof}

\medskip We can now prove Theorem~\ref{A} from the Introduction:

\begin{Theorem} \label{hyperbolic}
  Let $\mathbb A$ be a finite graph of groups where each vertex group
  either is polycyclic-by-finite or is word-hyperbolic and locally
  quasiconvex, and where all edge groups are virtually polycyclic.
  Then for any $v_0\in VA$ the group $G=\pi_1({\mathbb A}, v_0)$ has
  solvable uniform membership problem.  Moreover there is an algorithm
  which, given a finite subset $S\subseteq G$, constructs the induced
  splitting and a finite presentation for the subgroup $U=\langle
  S\rangle\le G$.
\end{Theorem}

\begin{proof}

  The uniform membership problem is solvable in polycyclic-by-finite
  groups~\cite{BCRD} and in locally quasiconvex hyperbolic
  groups~\cite{Ka96}.  Thus by Theorem~\ref{member} to establish the
  solvability of the membership problem it suffices to check that the
  graph of groups $\mathbb A$ is benign.

  It is well-known (see for example \cite{CDP90}) that a polycyclic
  subgroup of a word-hyperbolic group is virtually cyclic. Hence all
  edge groups for edges incident to hyperbolic vertex groups are in
  fact virtually cyclic.

  Suppose first that $v\in VA$ is such that $A_v$ is word-hyperbolic
  and locally quasiconvex.  Let $L$ be the regular language of all
  Short-Lex geodesic words in $A_v$ over some fixed finite generating
  set of $A_v$. It is well known that $L$ gives a bi-automatic
  structure with uniqueness for $A_v$. Since $A_v$ is assumed to be
  locally quasiconvex, all finitely generated subgroups of $A_v$ are
  $L$-rational.  Therefore by the result of \cite{Ka96}, there is a
  uniform algorithm which, given a finite set $X\subseteq A_v$,
  produces the pre-image $L_X$ of the subgroup $\langle X\rangle\le
  A_v$ in $L$. For each edge $e\in EA$ with $o(e)=v$ denote by $L_e$
  the pre-image in $L$ of the virtually cyclic subgroup
  $\alpha_e(A_e)$.

  Suppose now that $X\subseteq A_v$ is a finite set, $a\in A_v$ and
  $e\in EA$ is an edge with $o(e)=v$. We first construct the language
  $L_X$. Then using the biautomatic structure on $A_v$ we construct
  the regular language $L_{e,a}$ which is the pre-image in $L$ of the
  set $a\alpha_e(A_e)$. Now to decide if $\langle X\rangle \cap
  a\alpha_e(A_e)$ is empty we only need to check whether the
  intersection of the regular languages $L_X \cap L_{e,a}$ is empty.

  Moreover, we can also compute the intersection $L_X\cap L_e$ which
  is the pre-image in $L$ of the subgroup $\langle X\rangle \cap
  \alpha_e(A_e)$. Once the regular language $L_X\cap L_e$ is known, it
  is easy to recover a finite generating set for $\langle X\rangle
  \cap \alpha_e(A_e)$. Thus we have verified that $\mathbb A$ is
  benign at the vertex $v$.

  Suppose now that $A_v$ is virtually polycyclic. All virtually
  polycyclic groups are Noetherian and have solvable uniform
  membership problem (see for example \cite{BCRD}). Note that if
  $H,K\le A_v$ and $a\in A_v$ then $aH\cap K\ne \emptyset \iff a\in
  KH$. Since $A_v$ is virtually polycyclic, by a result of \cite{LW}
  the set $KH\subseteq A_v$ is closed in the profinite topology.
  Hence, given $a\in A_v$ and finite generating sets for $H,K$, we can
  detect if $a\not\in KH$ in some finite quotient of $A_v$. On the
  other hand, we can enumerate the set $KH$ and using the solvability
  of the word-problem in $A_v$, we can detect if $a\in KH$. Running
  this procedure parallel to enumerating all finite quotients of
  $A_v$, we can therefore decide whether or not $a$ belongs to $KH$.
  This shows that condition (1) of Definition~\ref{benign} holds at
  $v$. As proved in \cite{BCRD}, there is an algorithm which, given
  two finitely generated subgroups of a virtually polycyclic group,
  computes the generating set of their intersection.  Thus condition
  (4) of Definition~\ref{benign} also holds at $v$.

  We have verified that the graph of groups $\mathbb A$ is benign.
  Hence Theorem~\ref{member} applies and $G$ has solvable uniform
  membership problem.

  By Theorem~\ref{construct}, given a finite subset $S\subseteq G$ we
  can algorithmically construct a finite graph of group $\mathbb B$
  providing an induced splitting for $U=\langle S\rangle\le
  G=\pi_1({\mathbb A}, v_0)$. The vertex (and edge) groups of $\mathbb
  B$ are given as subgroups of vertex groups of $\mathbb A$ generated
  by some finite generating sets. By the result of
  Kapovich~\cite{Ka96} if $A_v$ is word-hyperbolic and locally
  quasiconvex, then there is an algorithm which, given a finite subset
  of $A_v$, produces a finite presentation for the subgroup generated
  by this set. The same is true for virtually polycyclic groups $A_v$,
  as proved in \cite{BCRD}. Hence we can recover a finite presentation
  of each vertex group of $\mathbb B$ and thus produce a finite
  presentation of $U$, as claimed.
\end{proof}

Not surprisingly, we also recover (a generalization of) Mihailova's
theorem regarding the membership problem for free products:

\begin{cor}\label{Mih}
  Let $G=\pi_1(\mathbb A,v_0)$, where $\mathbb A$ is a finite graph of
  finitely generated groups such that all edge groups are finite and
  all vertex groups have solvable membership problem. Then $G$ has
  solvable membership problem.

  Suppose, in addition, that for each $v\in VA$ there is an algorithm
  which, given a finite subset $Y$ of $A_v$, produces a finite
  presentation for the subgroup of $A_v$ generated by $Y$ (thus
  each $A_v$ is coherent). Then there is an algorithm which, given a
  finite set $S\subseteq G$, constructs a finite presentation for the
  subgroup $U=\langle S\rangle \le G$.
\end{cor}

\begin{proof}
  It is easy to see that $\mathbb A$ is benign and hence
  Corollary~\ref{Mih} follows from Theorem~\ref{construct} and
  Theorem~\ref{member}.
\end{proof}

\section{Grushko's Theorem}

As an application of our methods we can produce a quick proof of
Grushko's Theorem~\cite{Gru}. Recall that for a finitely generated
group $G$ \emph{the rank} $rk(G)$ is defined as the smallest number of
elements in a generating set of $G$.  A classical result of Grushko
states that rank behaves additively with respect to free products.

\begin{defn}[Complexity of an $\mathbb A$-graph]
  Let $\mathcal B$ be a finite $\mathbb A$-graph.

  We define \emph{the complexity of} $\mathcal B$ as
\[
c({\mathcal B}):=rk(\pi_1(B))+ \sum_{u\in VB} rk (B_u).
\]
\end{defn}
Recall that $\pi_1(B)$ is a free group whose rank is equal to the
number of edges in the complement of any maximal subtree of $B$.

\begin{Theorem}[Grushko~\cite{Gru}]
  Let $G_1, G_2$ be nontrivial finitely generated groups. Then
\[
rk(G_1 \ast G_2) =rk(G_1)+rk(G_2)
\]
\end{Theorem}
\begin{proof}
  It is obvious that $rk(G_1 \ast G_2) \le rk(G_1)+rk(G_2)$. Thus it
  suffices to establish the opposite inequality.

  Consider an edge of groups $\mathbb A$ with a single edge $e$, two
  vertices $v_0=o(e), v_1=t(v)$, the trivial edge group $A_e=1$ and
  vertex groups $A_{v_1}=G_1$ and $A_{v_2}=G_2$. Then
  \[G:=\pi_1({\mathbb A}, v_0)=G_1\ast G_2.\]

  Let $S$ be a generating set of $G$ of minimal cardinality, given as
  a collection of $\mathbb A$-reduced paths from $v_0$ to $v_0$. Thus
  $\# S=rk(G)$. Put $({\mathcal B_0}, u_0)$ to be the $S$-wedge. Notice
  that by construction $c({\mathcal B_0})=\# S$. We then start the
  abstract folding algorithm and construct a sequence of $\mathbb
  A$-graphs $({\mathcal B_0}, u_0), ({\mathcal B_1}, u_1), \dots $ by
  performing folding moves. Since the edge group in $\mathbb A$ is
  trivial, moves of type $F5-F6$ will never occur. Each move of type
  $F1-F4$ reduces the number of edges, and hence this sequence will
  terminate with a folded graph $({\mathcal B_n}, u_n)$. It is easy to
  see that moves $F1-F4$ do not increase the complexity and so
  $c({\mathcal B_n})\le c({\mathcal B_0})=\#S=rk(G)$. On the other
  hand $\mathcal B_n$ provides the induced splitting for the subgroup
  generated by $S$, that is for $G$ itself. Thus $\mathcal B_n$ recovers the original splitting $\mathbb A$ of $G$ which implies that $c({\mathbb B}_n)=rk (G_1)+rk(G_2)$. Thus $rk (G_1)+rk(G_2)\le rk(G)$.\end{proof}

\end{document}